\documentclass[SIADSpaper,onefignum,onetabnum]{siamonline171218}

\usepackage{lipsum}
\usepackage{amsfonts}
\usepackage{graphicx}
\usepackage{epstopdf}
\usepackage{algorithmic}
\ifpdf
  \DeclareGraphicsExtensions{.eps,.pdf,.png,.jpg}
\else
  \DeclareGraphicsExtensions{.eps}
\fi

\usepackage{bm}
\usepackage{color}
\usepackage{amsmath}
\usepackage{amssymb}
\usepackage{cleveref}
\usepackage{amssymb}
\usepackage[USenglish]{babel}
\usepackage{float}
\usepackage{graphicx}
\usepackage{mathtools}
\graphicspath{ {images/} }

\DeclareMathOperator*{\argmin}{\arg\!\min}

\usepackage{enumitem}
\setlist[enumerate]{leftmargin=.5in}
\setlist[itemize]{leftmargin=.5in}

\newsiamremark{remark}{Remark}
\newsiamremark{hypothesis}{Hypothesis}
\crefname{hypothesis}{Hypothesis}{Hypotheses}
\newsiamthm{claim}{Claim}
\usepackage{pifont}
\newcommand{\cmark}{\ding{51}}%
\newcommand{\xmark}{\ding{55}}%

\headers{Physics-Informed Probabilistic Learning of Linear Embeddings of Non-linear Dynamics With Guaranteed Stability}{S. Pan, and K. Duraisamy}

\title{Physics-Informed Probabilistic Learning of Linear Embeddings of Non-linear Dynamics With Guaranteed Stability\thanks{Submitted to the editors DATE.
\funding{DARPA Physics of AI Program}}}

\author{
Shaowu Pan\thanks{Department of Aerospace Engineering, University of Michigan, Ann Arbor, MI (\email{shawnpan@umich.edu}, \email{kdur@umich.edu} )}
\and
Karthik Duraisamy\footnotemark[2]
}

\usepackage{amsopn}

\ifpdf
\hypersetup{
  pdftitle={Physics-Informed Probabilistic Learning of Linear Embeddings of Non-linear Dynamics With Guaranteed Stability},
  pdfauthor={Shaowu Pan and Karthik Duraisamy}
}
\fi

\begin{document}

\maketitle

\begin{abstract}

The Koopman operator has emerged as a powerful tool for the analysis of nonlinear dynamical systems as it provides coordinate transformations to  globally linearize the dynamics. While recent deep learning approaches have been useful in extracting the Koopman operator from a data-driven perspective, several challenges remain.
In this work, we formalize the problem of learning the continuous-time Koopman operator with deep neural networks in a measure-theoretic framework. Our approach induces two types of models: differential and recurrent form, the choice of which depends on the availability of the governing equations and data. We then enforce  a structural parameterization that renders the realization of the Koopman operator provably stable. A new autoencoder architecture is constructed, such that only the residual of the dynamic mode decomposition is learned. Finally,  we employ mean-field variational inference (MFVI) on the aforementioned framework in a hierarchical Bayesian setting to quantify  uncertainties in the characterization and prediction of the dynamics  of observables.
The framework is evaluated on a simple polynomial system, the  Duffing oscillator, and an unstable cylinder wake flow with noisy measurements. 


\end{abstract}

\begin{keywords}
  Koopman decomposition, autoencoders, linear embedding,  reduced-order modeling, variational inference
\end{keywords} 

\begin{AMS}
  37M10, 37M25, 47B33, 62F15
\end{AMS}

%


\section{Introduction}

Nonlinear dynamical systems are prevalent in physical, mathematical and engineering problems~\cite{wiggins2003introduction}. Nonlinearity gives rise to rich phenomena, and is much more challenging to characterize in  comparison to linear systems for which comprehensive techniques have been developed for decades~\cite{khalil1996noninear}. 

Consider the autonomous system $\dot{\mathbf{x}} = \mathbf{F}(\mathbf{x})$, $\mathbf{x} \in \mathcal{M} \subset \mathbb{R}^N$,  where  $\mathcal{M}$ is a smooth manifold in the state space, $\mathbf{F}: \mathcal{M} \mapsto T\mathcal{M}$ is a vector-valued smooth function and $T\mathcal{M}$ is the tangent bundle, i.e., $\forall p \in \mathcal{M}, \mathbf{F}(p) \in T_p \mathcal{M} $. Instead of a geometric viewpoint~\cite{guckenheimer2013nonlinear,hirsch2012differential}, Koopman~\cite{koopman1931hamiltonian} offered an operator-theoretic perspective by describing the ``dynamics of observables" on the measure space $(\mathcal{M}, \Sigma, \mu)$ via the Koopman operator $\mathcal{K}_t: \mathcal{F} \mapsto \mathcal{F}$, such that for an observable on the manifold, $\forall f \in \mathcal{F}$, $t \in \mathbb{R}_{\ge 0}$, $f: \mathcal{M} \mapsto \mathbb{C}$, $\mathcal{K}_t f \triangleq f \circ S_t $, where $S_t(\mathbf{x}_0)\triangleq S(\mathbf{x}_0, t)$ is the flow map that takes the initial condition $x_0$ and advances it by time $t$ by solving the initial value problem for the aforementioned nonlinear dynamics, with $\mathcal{F} = L^2(\mathcal{M},\mu)$, where $\mu$ is some positive measure on $(\mathcal{M}, \Sigma)$. 

The semigroup of Koopman operators $\{\mathcal{K}_t\}_{t\in\mathbb{R}^+}$ is generated by  $\mathcal{K}: \mathcal{D}(\mathcal{K}) \mapsto \mathcal{F}$, $\mathcal{K}f \triangleq \lim_{t\rightarrow 0} (\mathcal{K}_t f - f)/t$  where $\mathcal{D}(\mathcal{K})$  is the domain in which the aforementioned limit is well-defined and $\overline{\mathcal{D}(\mathcal{K})} = \mathcal{F}$. This operator is {\em linear}, but defined on $\mathcal{F}$, which makes it inherently {\em infinite-dimensional}. 
Physically, $\mathcal{K}_t$ governs the temporal evolution of any observable in $\mathcal{F}$.  

Specifically, for  measure-preserving systems, e.g., Hamiltonian system or the dynamics on an attractor of the Naiver--Stokes equations, one can guarantee a spectral decomposition of any observable $f \in \mathcal{F}$, or the existence of a spectral resolution of the Koopman operator restricted on the attractor, i.e., $\mathcal{K}_t|_{\mathcal{A}}$, where $\mathcal{A} \subset \mathcal{M}$ is the attractor. Such a spectral decomposition~\cite{mezic2005spectral} is a sum of three essentials: temporal averaging of the observable, the contribution from the point spectrum, which is \emph{almost periodic} in time, and that from the continuous spectrum, which is \emph{chaotic}~\cite{koopman1932dynamical}.  For a more comprehensive discussion, readers are referred to the excellent review by Budi\v{s}i\'{c} et al.~\cite{budivsic2012applied}. It should be stressed that for measure-preserving systems, the Koopman operator is not only well-defined on $L^2(\mathcal{M},\mu)$ but also can be shown to be unitary, which ensures  properties such as simple eigenvalues and the existence of the aforementioned spectral resolution~\cite{budivsic2012applied}. Throughout this work,  we employ the simple but practical assumption~\cite{korda2018convergence} $\mathcal{F} = L^2(\mathcal{M}, \mu)$, where $\mu$ is some positive measure with support equal to $\mathcal{M}$. Note that this implies that the Koopman operator is well-defined on $\mathcal{F}$. 
 
The appeal of the Koopman operator lies in the possibility of globally linearizing the nonlinear dynamics of certain observables, i.e., the quantities of interest. In practice, we are mostly interested in the subspace $\mathcal{F}_D$, which consists of a finite number of Koopman eigenfunctions, namely the \emph{finite-dimensional Koopman invariant subspace}, where $D \in \mathbb{N}$ is the dimension of the subspace. Naturally, for a given observable $f$, we are  interested in the minimal $\mathcal{F}_D$, such that $f \in \mathcal{F}_D$. Koopman analysis entails the search for the Koopman eigenfunctions that span $\mathcal{F}_D$ and its associated eigenvalues. While we assume that the eigenvalues are simple, it is possible to extend to the generalized case~\cite{budivsic2012applied}. The identity mapping, i.e., $f = \mathbf{x}$ is of central interest since it corresponds to a linear system that is topologically conjugate to the nonlinear dynamical system. This essentially generalizes the Hartman-Grobman theorem to the entire basin of attraction of the equilibrium point or periodic orbit~\cite{lan2013linearization}.  From the viewpoint of understanding the behavior of the dynamical system in state space, the level sets of Koopman eigenfunctions form the invariant partition in the state space~\cite{budivsic2012geometry}, which can help study mixing. Further, the ergodic partition can be analyzed with Koopman eigenfunctions~\cite{mezic1999method}.  

 Koopman analysis has also gained increased attention in different communities.  Koopman eigenfunctions have been leveraged in nonlinear optimal control~\cite{kaiser2017data,mezic2015applications,brunton2016koopman,proctor2018generalizing,korda2016linear}. In modal analysis, the equivalence between the dynamic mode decomposition and the approximation of  Koopman modes was recognized by Rowley et al.~\cite{rowley2009spectral}. In a reduced order modeling context, the operator-theoretic viewpoint of the Koopman operator is critical to the formulation of Mori-Zwanzig (M-Z) formalism~\cite{chorin}, which has recently found utility in addressing closure and stabilization~\cite{gouasmi2017priori,parish2017dynamic, parish2017non} of multiscale problems. 






 Despite its appealing properties, the  Koopman decomposition cannot be pursued in its original form described above for practical scientific applications, because the operator is defined on the infinite dimensional Hilbert space. To accommodate practical computation of the  Koopman operator, we consider two simplifications:

\begin{itemize}

\item \emph{Only a finite dimensional space invariant to  $\mathcal{K}_t$ is of interest.} This excludes the possibility of dealing with a chaotic system since it is impossible for a finite dimensional linear system to be topologically conjugate to a chaotic system~\cite{budivsic2012applied}. However, one could assign a state-dependent Koopman operator to work around this limitation~\cite{lusch2017deep}.

\item \emph{Nonlinear reconstruction of  observables from Koopman eigenfunctions.}
The original concept of the Koopman operator studies the evolution of any observable in the Koopman invariant subspace. Thus, the observables of interest can be linearly reconstructed from the Koopman eigenfunctions. 
Regarding the latter Linear reconstruction is desirable especially in the content of control~\cite{kaiser2017data}. Li et al.~\cite{li2017extended} considered augmenting the Koopman invariant subspace with neural network-trained observables together with the system state $\mathbf{x}$ to force linear reconstruction.
However, for systems with multiple fixed points, there does not exist a finite-dimensional Koopman invariant subspace that also \emph{spans} the state \emph{globally}\footnote{\textcolor{black}{In this  paper, the term ``global embedding" is used to imply non-locality in phase space. Note that the Hartman-Grobman theorem~\cite{arrowsmith1992dynamical} establishes  a topological conjugacy to a linear system with the same eigenvalues in a small neighborhood of the fixed point.} } due to the impossibility of establishing such a topological conjugacy~\cite{brunton2016koopman} with the linear system. Therefore, attention naturally moves to the introduction of nonlinear reconstruction that is more expressive than linear reconstruction, together with extra modes that indicate different basins of each attractor~\cite{takeishi2017learning,otto2017linearly}. Unfortunately, nonlinear reconstruction is equivalent to removing Koopman modes for the observables of interest~\cite{otto2017linearly}. In other words, this is equivalent to removing the constraint that the system state lies in the finite dimensional Koopman invariant subspace. Although for systems with a single attractor, one could still retain a linear reconstruction representation, which \emph{might} lead to a decrease in accuracy~\cite{otto2017linearly}. 
\end{itemize}

To determine the Koopman eigenfunctions, we consider addressing the problem from a data-driven perspective. 
A common approximation is to assume that the Koopman invariant subspace $\mathcal{F}_D$ is spanned by a finite dictionary of  functions such as monomials, and then minimize the empirical risk of the residual that comes from the imperfect $\mathcal{F}_D$. This is the general idea behind Extended Dynamic Mode Decomposition (EDMD)~\cite{williams2015data} and Kernel-based Dynamic Mode Decomposition (KDMD)~\cite{williams2014kernel} where an \emph{implicitly}-defined infinite dimensional Koopman invariant subspace can be computationally affordable if the kernel satisfies Mercer's condition~\cite{j1909xvi}. Note that  EDMD is essentially a least-squares regression of the (linear) action of the Koopman operator for the features in the dictionary, with samples drawn from some measure $\mu$ with fixed features, which can be proved as a $L^2(\mathcal{M},\mu)$ orthogonal projection in the limit of large independent identically distributed data~\cite{klus2015numerical,williams2015data}. 

By framing the data-derived Koopman operator in the Hilbert space endowed with a proper measure, one can prove  optimality in the asymptotic sense for EDMD~\cite{korda2018convergence} and Hankel-DMD~\cite{arbabi2017ergodic}. As a side note, if the dynamics on the attractor is of central interest, it has been shown that Hankel-DMD with appropriate delay embedding and sampling rate can recover the entire system, even without sampling the full attractor~\cite{pan2019structure}. Further, one can establish a duality between the measure-preserving dynamics on the attractor with a stationary stochastic process, which reflects a close relationship between the estimation of a continuous spectrum from the trajectories, and the extraction of the power spectral density of stochastic signals~\cite{arbabi2017study}. Other algorithms include generalized Laplace analysis (GLA)~\cite{budivsic2012applied,mauroy2012use}, and the Ulam Galerkin method~\cite{froyland2014computational}. If the representation of the Koopman eigenfunctions is \emph{sparse} in the pre-defined dictionary in EDMD, then one can leverage sparse regression techniques~\cite{brunton2016discovering} to select observables in an iterative fashion~\cite{brunton2016koopman}. 
In this work, we do not address continuous spectra. One way of addressing continuous spectra is to use~\cite{lusch2017deep} an auxiliary sub-network to characterize the eigenvalues as \emph{state-dependent}.  


\textcolor{black}{However, it has been shown that traditional EDMD/KDMD method is prone to overfitting and is hard to interpret~\cite{otto2017linearly}. For instance, the number of approximated Koopman eigenfunctions scales with either the number of features in the dictionary or the number of training snapshots even though only a few of them are genuine Koopman eigenfunctions.}
Recently, there have been several \textcolor{black}{alternative} attempts to leverage deep learning architectures~\cite{otto2017linearly,li2017extended,takeishi2017learning,wehmeyer2018time,yeung2017learning}, to extract  Koopman decompositions. Yeung et al.~\cite{yeung2017learning} used feedforward neural networks to learn the dictionary, but the reconstruction loss was found to be non-optimal. Li et al.~\cite{li2017extended} enforced several non-trainable functions, e.g., components of the system state, in the Koopman observables to ensure an accurate reconstruction but one that could be inefficient in terms of obtaining a finite-dimensional Koopman observable subspace~\cite{otto2017linearly}. Further, Takeishi et al.~\cite{takeishi2017learning} utilized  linear time delay embedding in the feedforward neural network framework to construct Koopman observables  with nonlinear reconstruction which is critical for partially observed systems~\cite{pan2018data}.  Lusch et al.~\cite{lusch2017deep} further extended the deep learning framework to chaotic systems. Otto and Rowley~\cite{otto2017linearly} considered a recurrent loss for better performance on long time prediction on trajectories that transit to the attractor. Recently, Morton et al.~\cite{morton2019deep} addressed the uncertainty in a deep learning model with a focus on control. \textcolor{black}{The benefit of formulating the search for the Koopman operator in an optimization setting enables the enforcement of stability. For example, it is also feasible to constrain eigenvalues in optimized DMD~\cite{chen2012variants,askham2018variable}. Specifically in the neural network context,} Erichson et al.~\cite{erichson2019physics} considered a stability promoting loss to encourage Lyapunov stability of the dynamical system. 

A unified approach towards uncertainty quantification, stabilization, and incorporation of physics information is lacking, particularly in a continuous setting. 
This motivates us to establish a probabilistic stabilized deep learning framework specifically to learn the Koopman decomposition for a continuous dynamical system. We employ automatic differentiation variational inference (ADVI)~\cite{kucukelbir2017automatic} to quantify parametric uncertainty in deep neural networks, and structural parameterization to enforce stability of the Koopman operator extracted from the deep neural networks. A broad comparison of the present work with other approaches is shown in \cref{tab:compare}.

\begin{table}[htbp]
{\footnotesize
\caption{Comparison with other approaches in the literature}
\label{tab:compare}
\begin{center}
 \begin{tabular}{c c c c c c } 
 Previous works & \shortstack{continuous\\/discrete} & \shortstack{nonlinear\\reconstruction}& \shortstack{continuous\\spectrum} & uncertainty & stability \\ [0.5ex] 
 \hline\hline
 Yeung et al.~\cite{yeung2017learning}    & discrete   & \textcolor{black}{\xmark} & \textcolor{black}{\xmark}  & \textcolor{black}{\xmark} & \textcolor{black}{\xmark}  \\ 
 Li et al.~\cite{li2017extended}       & discrete   &  \textcolor{black}{\xmark}  & \textcolor{black}{\xmark}  &  \textcolor{black}{\xmark} &  \textcolor{black}{\xmark} \\
 Takeishi et al.~\cite{takeishi2017learning}       & discrete   &  \textcolor{black}{\cmark}  & \textcolor{black}{\xmark}  &  \textcolor{black}{\xmark} &  \textcolor{black}{\xmark} \\
 Otto and Rowley~\cite{otto2017linearly}  & discrete   &  \textcolor{black}{\cmark}   & \textcolor{black}{\xmark} &  \textcolor{black}{\xmark} &  \textcolor{black}{\xmark} \\
 Lusch et al.~\cite{lusch2017deep}    & discrete   &  \textcolor{black}{\cmark} &  \textcolor{black}{\cmark}  &  \textcolor{black}{\xmark} &  \textcolor{black}{\xmark}  \\
  Morton et al.~\cite{morton2019deep}    & discrete   &  \textcolor{black}{\cmark} &  \textcolor{black}{\xmark}  &  \textcolor{black}{\cmark} &  \textcolor{black}{\xmark}  \\
 Erichson et al.~\cite{erichson2019physics}    & discrete   &  \textcolor{black}{\cmark} &  \textcolor{black}{\xmark}  &  \textcolor{black}{\xmark} &  \textcolor{black}{\cmark}  \\
\textcolor{black}{Our framework}          & \textcolor{black}{continuous} & \textcolor{black}{\cmark} & \textcolor{black}{\xmark}  & \textcolor{black}{\cmark} & \textcolor{black}{\cmark} \\ [1ex] 
 \end{tabular}
\end{center}
}
\end{table}   

The outline of this paper is as follows: In \cref{sec:cont_frame}, our  framework of deep learning of Koopman operators of continuous dynamics, together with structured parameterization that enforces stabilization is presented. In \cref{sec:bayesdl}, Bayesian deep learning, specifically  the employment of variational inference is introduced. Results and discussion on several problems are given in \cref{sec:results}. Conclusions are drawn in \cref{sec:conclusion}.

\section{Data-driven framework to learn continuous-time Koopman decompositions}
\label{sec:cont_frame}

 In contrast to  previous  approaches~\cite{otto2017linearly,lusch2017deep,yeung2017learning,li2017extended}, our framework pursues  continuous-time Koopman decompositions. The continuous formulation is more amenable to posit desired constraints and contend with non-uniform sampling~\cite{otto2017linearly,askham2018variable}, which is frequently encountered in experiments or temporal multi-scale data. To begin with, recall the general form of autonomous continuous nonlinear dynamical systems, 
\begin{equation}{\label{eq:nldy}}
\dot{\mathbf{x}} = \mathbf{F}(\mathbf{x}), \quad \mathbf{x} \in \mathcal{M} \subset \mathbb{R}^N.
\end{equation}
We seek a finite dimensional Koopman invariant subspace $\mathcal{F}_D$ with $D $ linearly independent smooth observation functions, defined as
\begin{equation}
\mathcal{F}_D = \textrm{span}\{ \phi_1, \ldots, \phi_D \} \subset \mathcal{F},
\end{equation}
where $\phi_i \in C^1(\mathcal{M}, \mathbb{R})$ and $\phi_i \in \mathcal{F}$, $i=1,\ldots, D$. 
Correspondingly, the observation vector $\mathbf{\Phi}$ is defined as
\begin{equation}{\label{eq:nlob}}
\mathbf{\Phi} (x) = \begin{bmatrix}
\phi_1(x) & \phi_2(x) & \ldots & \phi_D(x)
\end{bmatrix} \in \mathcal{F}^D.
\end{equation}
Based on the aforementioned structure, we require the following two conditions. First, $\mathbf{\Phi} (x) $ evolves linearly in time, i.e., $\exists$  $\mathbf{K} \in \mathbb{R}^{D \times D}$ s.t. 
\begin{equation}{\label{eq:lindy}}
\dot{\mathbf{\Phi}} \triangleq d\mathbf{\Phi}/dt =  \mathbf{\Phi} \mathbf{K}.
\end{equation} 
By the chain rule, the relationship between \cref{eq:lindy} and \cref{eq:nldy} is
\begin{equation}{\label{eq:keyeq}}
\mathbf{F} \cdot \nabla_{\mathbf{x}} {\mathbf{\Phi}} = \mathbf{\Phi} \mathbf{K}.
\end{equation}
Second, $\exists$ $\mathbf{\Psi}: \mathbb{R}^D \mapsto \mathcal{M}$, s.t. $ \mathbf{\Psi} \circ \mathbf{\Phi} = \mathcal{I}$, $\mathcal{I}: \mathcal{M} \mapsto \mathcal{M}$ is the identity map. Therefore, we can recover the state $\mathbf{x}$ from $\mathbf{\Phi}$. 

The goal is to find a minimal set of basis functions that  spans a Koopman invariant subspace. There are several  existing methods for approximation: 
\begin{itemize}
\item Dynamic mode decomposition (DMD): observation functions are linear transformations of the state, usually POD modes for robustness purposes~\cite{schmid2010dynamic}. 
\item Extended/Kernel DMD: observation functions are pre-specified, either explicit polynomials~\cite{williams2015data} or implicitly by the kernel~\cite{williams2014kernel}. 
\item Deep learning Koopman/Neural Network DMD~\cite{otto2017linearly,lusch2017deep} searching a set of fixed size for the nonlinear observations by artificial neural networks.
\end{itemize}

In this work, we specifically focus on using neural networks. 
\subsection{Neural network}

 The basic building block we use is  a feedforward neural network (FNN). A typical FNN with $L$ layers is the mapping $g(\cdot; \mathbf{W}_g): \mathbb{R}^n \mapsto \mathbb{R}^m$ such that
\begin{equation}
\eta_l = \sigma_l (\eta_{l-1} W_l  + b_l),
\end{equation}
for $l=1,\dots,L-1$, where $\eta_{0}$ stands for the input of the neural network $x$, $\eta_l \in \mathbb{R}^{n_l}$ is the number of hidden units in the layer $l$, $\sigma_l$ is activation function of layer $l$. Note that the last layer is linear, i.e., $\sigma_L(x)=x$: 
\begin{equation}
g( \eta_0;  \mathbf{W}_g)  = \eta_{L-1} W_L   + b_L, 
\end{equation}
where $\mathbf{W}_g  = \{W_1, b_1,\ldots, W_L, b_L\}$, $W_l \in \mathbb{R}^{n_{l-1} \times n_{l}}$, $b_l \in \mathbb{R}^{n_l}$, for $l=1, \ldots, L$. Note that $n_0 = n \in \mathbb{N}$, $n_L=m  \in \mathbb{N}$. Such a mapping, given an arbitrary number of hidden units, even with a single hidden layer, has been shown to be a universal approximator in the $L^p$ sense $1 \le p < \infty$, as long as the activation function is not polynomial almost everywhere~\cite{leshno1993multilayer}. Throughout this work, we use the Swish~\cite{ramachandran2017searching} activation function, which is continuously differentiable and found to achieve strong empirical performance over many variants of ReLU and ELU on typical deep learning tasks.

\subsection{Problem formulation}

We first define the ideal problem of learning the Koopman decomposition  as a constrained variational problem, and incorporate assumptions to make it tractable step by step. Then, we introduce the effect of data as an empirical measure into the optimization in the function space. 
Specifically, we propose two slightly different formulations based on varying requirements. 
\begin{enumerate}
\item Differential form: for low-dimensional systems when the governing equations, i.e., \cref{eq:keyeq} are known, they can be leveraged  without the trajectory data.
\item Recurrent form: for high-dimensional systems where only discrete trajectory data is obtained, which implies the access to the action of $\mathcal{K}_t$ over discrete $t$. 
\end{enumerate}

Recall that we are searching the Koopman operator defined on  $\mathcal{F} = L^2(\mathcal{M},\mu)$, which is the space of all measurable functions $\phi: \mathcal{M} \mapsto \mathbb{R}$ such that,  
\begin{equation}
\left\lVert \phi \right\rVert_{\mathcal{F}} \triangleq \sqrt{\int_{\mathcal{M}} |\phi|^2 d\mu} < \infty.
\end{equation} 
As a natural extension, for any finite $ n$, given the vector-valued function $\mathbf{\Phi} = \begin{bmatrix}
\phi_1 & \ldots & \phi_n
\end{bmatrix} \in \mathcal{F}^n$, we define the corresponding norm as,
\begin{equation}
\left\lVert \mathbf{\Phi} \right\rVert_{\mathcal{F}^n} \triangleq \sqrt{\int_{\mathcal{M}} \sum_{j=1}^{n} |\mathbf{\phi}_j|^2 d\mu}.
\end{equation}

\subsubsection{Leveraging known physics: differential form}
\label{sec:diff_form}

For any observation vector $\mathbf{\Phi} \in \mathcal{F}^D$, we can define the following Koopman error functional $\mathcal{J}[\cdot]: \mathcal{F}^D \mapsto \mathbb{R}_{\ge 0}$, $\mathbb{R}_{\ge 0} = [0,+\infty)$ as the square of maximal distance for all components in $\Phi$ between the action of $\mathcal{K}$ and its $L^2$ projection onto $\mathcal{F}_D$~\cite{korda2018convergence}, 
\begin{equation}
\mathcal{J}[\mathbf{\Phi}] = \max_{\psi \in \{ \phi_1, \ldots, \phi_D \}} \min_{f \in \mathcal{F}_D} \left\lVert f - \mathcal{K} \psi \right\rVert^2_{\mathcal{F}},
\end{equation}
which describes the extent to which $\mathcal{F}_D$ is invariant to the Koopman operator $\mathcal{K}_t$ with $t \rightarrow 0$ with respect to each basis. Ideally, if we can find $\mathbf{\Phi}$ such that the corresponding $\mathcal{F}$ is invariant to $\mathcal{K}$, i.e., $\mathcal{J}[\mathbf{\Phi}] = 0$, then $\mathcal{F}_D$ is invariant to $\mathcal{K}_t$, $\forall t >0$, i.e., a perfect \emph{linear} embedding in the $L^2(\mu)$ sense. Once such an embedding is determined, the realization of $\mathcal{K}$ is simply the matrix $\mathbf{K}$. In this work,  we are interested in  $\mathbf{\Phi}$ such that one can recover $\mathbf{x}$ and $\mathcal{J}[\mathbf{\Phi}] \ge 0$. 

Although the above problem setup only contains minimal necessary assumptions, it is both mathematically and computationally challenging. For practical purposes, following previous  studies~\cite{otto2017linearly,williams2014kernel,williams2015data,schmid2010dynamic,lusch2017deep,takeishi2017learning},  we consider the following  assumptions to make the problem tractable:
\begin{enumerate}

\item Instead of solving the equation $\mathcal{J}[\mathbf{\Phi}] = 0$, we search for $\mathbf{\Phi}$ by finding the minimum of the following constrained variational problem,
\begin{equation}
\mathbf{\Phi}^* = \argmin_{\substack{\mathbf{\Phi} \in \mathcal{F}^D, \exists \mathbf{\Psi}: \mathbb{R}^D \mapsto \mathcal{M} \\ \mathbf{\Psi} \circ \mathbf{\Phi} = \mathcal{I} }} \mathcal{J}[\mathbf{\Phi}].
\end{equation}

\item Instead of directly solving the variational problem in the infinite dimensional $\mathcal{F}^D$, we optimize $\mathbf{\Phi}$ in the finite dimensional parameter space of feedforward neural networks with fixed architecture in which the number of layers and the number of hidden units in each layer is determined heuristically. Note that we are searching in a subset of $\mathcal{F}^D$ described by $\mathbf{W}_{\mathbf{\Phi}}$, which might induce a gap due to the choice of the neural network architecture, 
\begin{equation}
0 \le \min_{\substack{\mathbf{\Phi} \in \mathcal{F}^D, \exists \mathbf{\Psi}: \mathbb{R}^D \mapsto \mathcal{M} \\ \mathbf{\Psi} \circ \mathbf{\Phi} = \mathcal{I} }} \mathcal{J}[\mathbf{\Phi}] \le 
\min_{\substack{\mathbf{W}_{\mathbf{\Phi}}, \exists \mathbf{\Psi} \in C(\mathbb{R}^D, \mathbb{R}^N) \\ \mathbf{\Psi} \circ \mathbf{\Phi}(\cdot; \mathbf{W}_{\mathbf{\Phi}}) = \mathcal{I} }} \mathcal{J}[\mathbf{\Phi}(\cdot; \mathbf{W}_{\mathbf{\Phi}})].
\end{equation}
Clearly, the gap is bounded by the right hand side of the second inequality above. In addition, it should be noted that the requirement of linear independence in $\{\phi_1,\ldots,\phi_D\}$ is relaxed, but $\dim \mathcal{F}_D $ is bounded by $D$.

\item As a standard procedure in deep learning~\cite{goodfellow2016deep}, we use the penalty method  to \emph{approximate} the constrained optimization problem with an unconstrained optimization with unity penalty coefficient. Since this still entails nonconvex optimization, we define a global minima  as follows:
\begin{align}
\mathbf{W}_{\mathbf{\Phi}}^* &= \argmin_{\substack{\mathbf{W}_{\mathbf{\Phi}}, \exists \mathbf{\Psi} \in C(\mathbb{R}^D, \mathbb{R}^N) \\ \mathbf{\Psi} \circ \mathbf{\Phi}(\cdot; \mathbf{W}_{\mathbf{\Phi}}) = \mathcal{I} }} \mathcal{J}[\mathbf{\Phi}(\cdot; \mathbf{W}_{\mathbf{\Phi}})],  \\
\widehat{\mathbf{W}}_{\mathbf{\Phi}}, \widehat{\mathbf{W}}_{\mathbf{\Psi}} &=  \argmin_{\substack{\mathbf{W}_{\mathbf{\Phi}}, \mathbf{W}_{\mathbf{\Psi}}}} \mathcal{J}[\mathbf{\Phi}(\cdot; \mathbf{W}_{\mathbf{\Phi}})] +  \mathcal{R}[\mathbf{\Phi}(\cdot; \mathbf{W}_{\mathbf{\Phi}}), \mathbf{\Psi}(\cdot; \mathbf{W}_{{\mathbf{\Psi}}}) ],
\end{align}
where the reconstruction error functional $\mathcal{R}[\cdot,\cdot]: \mathcal{F}^D \times C(\mathbb{R}^D, \mathbb{R}^N) \mapsto \mathbb{R}_{\ge 0}$, is defined as 
\begin{equation}
\mathcal{R}[\mathbf{\Phi}, \mathbf{\Psi}] = \left\lVert \mathbf{\Psi} \circ \mathbf{\Phi}   - \mathcal{I} \right\rVert^2_{\mathcal{F}^N},
\end{equation}
for $\mathbf{\Phi} \in \mathcal{F}^D$, $\mathbf{\Psi} \in C(\mathbb{R}^D, \mathbb{R}^N) $. Then we assume
$\mathbf{\Phi}(\cdot; \widehat{\mathbf{W}}_{\mathbf{\Phi}})$ approximates one of the global minima, i.e., $\mathbf{\Phi}(\cdot; \mathbf{W}^*_{\mathbf{\Phi}})$.
Note that the convergence to a global minimum of the constrained optimization can be guaranteed if one is given the global minima of the sequential unconstrained optimization and by increasing the penalty coefficient to infinity~\cite{luenberger1973introduction}. 

\item We then optimize the sum of square error $\widehat{\mathcal{J}}$ over all components of $\mathbf{\Phi}$, which serves a upper bound of $\mathcal{J}$, 
\begin{align}
\mathcal{J}[\mathbf{\Phi}] \le \widehat{\mathcal{J}}[\mathbf{\Phi}] &= { \sum_{i=1}^{D} \min_{f \in \mathcal{F}_D} \left\lVert f - \mathcal{K} \phi_i \right\rVert^2_{\mathcal{F}}} = \min_{\mathbf{K} \in \mathbb{R}^{D\times D}}  \widetilde{\mathcal{J}}[\mathbf{\Phi}, \mathbf{K}], \\
\widetilde{\mathcal{J}}[\mathbf{\Phi}, \mathbf{K}] &=  \left\lVert \mathbf{\Phi} \mathbf{K} - \mathcal{K} \mathbf{\Phi} \right\rVert^2_{\mathcal{F}^D}.
\end{align}
The above formulation also implies  equal importance among all components of $\mathbf{\Phi}$. 

\item Despite the non-convex nature of the problem, we employ gradient-descent optimization to search for a local minimum. 

\end{enumerate}

In summary based on the above assumptions, we will solve the following optimization problem by gradient-descent:
\begin{equation}
\label{eq:diff_objective}
\widehat{\mathbf{W}}_{\mathbf{\Phi}}, \widehat{\mathbf{W}}_{\mathbf{\Psi}}, \widehat{\mathbf{K}} =  \argmin_{\substack{\mathbf{W}_{\mathbf{\Phi}}, \mathbf{W}_{\mathbf{\Psi}}, \mathbf{K}}} \widetilde{\mathcal{J}}[\mathbf{\Phi}(\cdot; \mathbf{W}_{\mathbf{\Phi}}), \mathbf{K}] +  \mathcal{R}[\mathbf{\Phi}(\cdot; \mathbf{W}_{\mathbf{\Phi}}), \mathbf{\Psi}(\cdot; \mathbf{W}_{{\mathbf{\Psi}}}) ].
\end{equation}

\subsubsection{Unknown physics with only the trajectory data: recurrent form}
\label{sec:recurrent}

From the viewpoint of approximation, it can be argued that the most natural way to determine the continuous Koopman operator is in the aforementioned differential form. 
However, higher accuracy can be achieved by taking advantage of trajectory information and minimizing the error over multiple time steps in the trajectory. This is the key idea behind optimized DMD~\cite{askham2018variable}. Recently, Lusch et al.~\cite{lusch2017deep}, Otto and Rowley~\cite{otto2017linearly} extended this idea to determine the discrete-time Koopman operator using deep learning. 

Recall that we assume that the space of observation functions $\mathcal{F} = L^2(\mathcal{M}, \mu)$  is invariant to $\mathcal{K}_t$, $\forall t \in \mathbb{R}_{\ge 0}$. Thus, we consider the $t$-time evolution of any function $\phi \in \mathcal{F}$, i.e., $\mathcal{K}_t \phi(\mathbf{x})$, as a function on $\mathcal{U} = \mathcal{M} \times \mathcal{T}$,  given an initial condition $\mathbf{x} \in \mathcal{M}$ and time evolution $t \in \mathcal{T}$ where $\mathcal{T}$ is a one dimensional smooth manifold, sometimes also referred to as the time horizon. We assume $\mathcal{K}_t \phi(\mathbf{x}) \in \mathcal{G} = L^2(\mathcal{U}, \nu)$. Based on the fact that $\mathcal{F}$ is invariant to $\mathcal{K}_t$, such an assumption can be shown to be valid  for compact $\mathcal{T}$ for a proper measure $\nu$.



Similar to the differential form, for any observation vector $\mathbf{\Phi} \in \mathcal{F}^D$, we define the following Koopman error functional $\mathcal{J}_r[\cdot]: \mathcal{F}^D \mapsto \mathbb{R}_{\ge 0}$ as the square of maximal $L^2$ distance for all components in $\mathbf{\Phi}$ between the predicted and ground truth trajectory,  
\begin{equation}
\mathcal{J}_r[\mathbf{\Phi}] = \max_{\psi \in \{ \phi_1,\ldots, \phi_D\} } \min_{\mathbf{K} \in \mathbb{R}^{D \times D}} \left\lVert \mathbf{\Phi} e^{t\mathbf{K}} c_{\psi} - \mathcal{K}_t \psi \right\rVert^2_{\mathcal{G}},
\end{equation}
where $\psi=\mathbf{\Phi}c_{\psi}$, $c_{\psi} \in \mathbb{R}^{D \times 1}$.
Following similar assumptions in \cref{sec:diff_form},  we need to define a reconstruction error functional to describe the discrepancy between the reconstructed and original states. Indeed, one can directly define the prediction error functional in the recurrent form, $\widetilde{\mathcal{P}}[\cdot, \cdot, \cdot]: \mathcal{F}^D \times C(\mathbb{R}^D, \mathbb{R}^N) \times \mathbb{R}^{D\times D} \mapsto \mathbb{R}_{\ge 0}$ as
\begin{equation}
\widetilde{\mathcal{P}}[\mathbf{\Phi}, \mathbf{\Psi}, \mathbf{K}] = \left\lVert
\mathbf{\Psi} \circ \mathbf{\Phi}  e^{t \mathbf{K}} - \mathcal{K}_t \mathcal{I}
\right\rVert^2_{\mathcal{G}^N}.
\end{equation}
Similarly, we define
\begin{equation}
\widetilde{\mathcal{J}}_r[\mathbf{\Phi}(\cdot; \mathbf{W}_{\mathbf{\Phi}}), \mathbf{K}] =  \left\lVert \mathbf{\Phi} e^{t \mathbf{K}} - \mathcal{K}_t \mathbf{\Phi} \right\rVert^2_{\mathcal{G}^D},
\end{equation}
and solve the following optimization problem via a gradient-based algorithm:
\begin{equation}
\label{eq:recurrent_objective}
\widehat{\mathbf{W}}_{\mathbf{\Phi}}, \widehat{\mathbf{W}}_{\mathbf{\Psi}}, \widehat{\mathbf{K}} =  \argmin_{\substack{\mathbf{W}_{\mathbf{\Phi}}, \mathbf{W}_{\mathbf{\Psi}}, \mathbf{K}}} \widetilde{\mathcal{J}}_r[\mathbf{\Phi}(\cdot; \mathbf{W}_{\mathbf{\Phi}}), \mathbf{K}] +  \widehat{\mathcal{P}}[\mathbf{\Phi}(\cdot; \mathbf{W}_{\mathbf{\Phi}}), \mathbf{\Psi}(\cdot; \mathbf{W}_{{\mathbf{\Psi}}}), \mathbf{K}].
\end{equation}

As a side note, one might also define the following reconstruction functional similar to previous differential form that is independent of $\mathbf{K}$, 
\begin{equation}
\mathcal{R}[\mathbf{\Phi}, \mathbf{\Psi}] = \left\lVert \mathcal{K}_t
(\mathbf{\Psi} \circ \mathbf{\Phi} - \mathcal{I}) 
\right\rVert^2_{\mathcal{G}^N},
\end{equation}
which can bound the prediction error functional together with the  Koopman error functional by triangular inequality, 
\begin{align}
\left\lVert
\mathbf{\Psi} \circ \mathbf{\Phi} e^{t\mathbf{K}} - \mathcal{K}_t \mathcal{I}
\right\rVert_{\mathcal{G}^N} &\le 
\left\lVert 
\mathbf{\Psi} \circ \mathbf{\Phi} e^{t\mathbf{K}} -  \mathcal{K}_t \mathbf{\Psi} \circ \mathbf{\Phi}
\right\rVert_{\mathcal{G}^N} +
\left\lVert \mathcal{K}_t
(\mathbf{\Psi} \circ \mathbf{\Phi} - \mathcal{I}) 
\right\rVert_{\mathcal{G}^N},\\
& \le L_{\mathbf{\Psi}} \left\lVert 
 \mathbf{\Phi} e^{t\mathbf{K}} -  \mathcal{K}_t \mathbf{\Phi}
\right\rVert_{\mathcal{G}^D} + 
\left\lVert \mathcal{K}_t
(\mathbf{\Psi} \circ \mathbf{\Phi} - \mathcal{I}) 
\right\rVert_{\mathcal{G}^N},
\end{align}
where $L_{\mathbf{\Psi}}$ is the Lipschitz constant for $\mathbf{\Psi}$. However, in this work, we do not have control over  $L_{\mathbf{\Psi}}$, and thus we prefer to directly minimize the prediction error function as in previous studies~\cite{otto2017linearly,lusch2017deep}.

\subsection{Measures}

\subsubsection{Measure for data generated by sampling in  phase space}

We consider the situation where $M$ data points on $\mathcal{M}$, i.e., $\{\mathbf{x}_m\}_{m=1}^{M}$ are drawn independently from some measure $\mu$ on $\mathcal{M}$, e.g., uniform distribution. This induces the following empirical measure $\hat{\mu}_M$, 
\begin{equation}
\hat{\mu}_M = \frac{1}{M} \sum_{m=1}^M \delta_{\mathbf{x}_m},
\end{equation}
where $\delta_x$ is the Dirac measure for $x$. Note that $\hat{\mu}_M$ uniformly converges to $\mu$~\cite{vapnik2015uniform} as $M \rightarrow \infty$. Thus, one can rewrite the differential form in \cref{eq:diff_objective} as an empirical risk minimization~\cite{vapnik1992principles}, 
\begin{align}
\label{eq:finite_diff}
&\widehat{\mathbf{W}}_{\mathbf{\Phi}}, \widehat{\mathbf{W}}_{\mathbf{\Psi}}, \widehat{\mathbf{K}} = \lim_{M \rightarrow \infty} \widehat{\mathbf{W}}_{\mathbf{\Phi},M}, \widehat{\mathbf{W}}_{\mathbf{\Psi},M}, \widehat{\mathbf{K}}_{M} \\ 
\nonumber &=  \lim_{M \rightarrow \infty}  \argmin_{\substack{\mathbf{W}_{\mathbf{\Phi}}, \mathbf{W}_{\mathbf{\Psi}}, \mathbf{K}}} \widetilde{\mathcal{J}}_{M}[\mathbf{\Phi}(\cdot; \mathbf{W}_{\mathbf{\Phi}}), \mathbf{K}] +  \mathcal{R}_{M}[\mathbf{\Phi}(\cdot; \mathbf{W}_{\mathbf{\Phi}}), \mathbf{\Psi}(\cdot; \mathbf{W}_{{\mathbf{\Psi}}}) ] \\
\nonumber &=  \lim_{M \rightarrow \infty} \argmin_{\substack{\mathbf{W}_{\mathbf{\Phi}}, \mathbf{W}_{\mathbf{\Psi}}, \mathbf{K}}}  \left\lVert \mathbf{\Phi}(\cdot; \mathbf{W}_{\mathbf{\Phi}}) \mathbf{K} - \mathcal{K} \mathbf{\Phi} \right\rVert^2_{\widehat{\mathcal{F}}_M^D} + \left\lVert \mathbf{\Psi}(\cdot; \mathbf{W}_{{\mathbf{\Psi}}}) \circ \mathbf{\Phi}(\cdot; \mathbf{W}_{\mathbf{\Phi}})   - \mathcal{I} \right\rVert^2_{\widehat{\mathcal{F}}_M^N}\\
\nonumber &=  \lim_{M \rightarrow \infty}  \argmin_{\substack{\mathbf{W}_{\mathbf{\Phi}}, \mathbf{W}_{\mathbf{\Psi}}, \mathbf{K}}} \frac{1}{M}\sum_{m=1}^{M}\bigg(  \lVert \mathbf{\Phi}(\mathbf{x}_m; \mathbf{W}_{\mathbf{\Phi}}) \mathbf{K} - \mathbf{F} \cdot \nabla_{\mathbf{x}} \mathbf{\Phi} (\mathbf{x}_m) \rVert^2  \\
&  \nonumber  \qquad \qquad  \qquad  \qquad + \left\lVert \mathbf{\Psi} (\mathbf{\Phi}(\mathbf{x}_m; \mathbf{W}_{\mathbf{\Phi}}); \mathbf{W}_{{\mathbf{\Psi}}})   - \mathbf{x}_m \right\rVert^2 \bigg),
\end{align}
where $\widehat{\mathcal{F}}_M = L^2(\mathcal{M}, \hat{\mu}_M)$.

\subsubsection{Measure for trajectory data generated by solving the initial value problem}

In the general case, information content  in trajectory data resulting from the solution of the initial value problem is strongly dependent on the initial condition. For instance, when the initial condition is in a region of phase space with sharp changes, it would be sensible to use a high sampling rate. On the other hand, if the initial condition is near a fixed point attractor, one might prefer to stop collecting the data after the system arrives at the fixed point. Such a sampling pattern for the specific initial state can be summarized as a Markov kernel~\cite{klenke2013probability} $\kappa: \mathcal{M} \times \Sigma_{\mathcal{T}} \mapsto [0,1]$. For every fixed initial state $\mathbf{x} \in \mathcal{M}$, the map $\kappa (\mathcal{E}_{\mathcal{T}}, \mathbf{x})$ is a measure on $(\mathcal{T}, \Sigma_{\mathcal{T}})$, $\mathcal{E}_{\mathcal{T}} \in \Sigma_{\mathcal{T}}$, where $\Sigma_{\mathcal{T}}$ is the $\sigma$-algebra on $\mathcal{T}$. Then, there exists a unique measure $\nu$ on $\mathcal{U}$~\cite{klenke2013probability} such that, 
\begin{equation}
\label{eq:markov_kernel}
\nu(\mathcal{E}_{\mathcal{M}} \times \mathcal{E}_{\mathcal{T}}) = \int_{\mathcal{E}_{\mathcal{M}}} \kappa(\mathcal{E}_{\mathcal{T}}, \mathbf{x}) d\mu(\mathbf{x}), 
\end{equation}
for $\mathcal{E}_{\mathcal{T}} \in \Sigma_{\mathcal{T}}, \mathcal{E}_{\mathcal{M}} \in \Sigma_{\mathcal{M}}$.

Assume we are given $M$ trajectories, $\{ \{\mathbf{x}_{m,j}\}_{j=1}^{T_m} \}_{m=1}^{M}$, $\{ \{t_{m,j}\}_{j=1}^{T_m} \}_{m=1}^{M}$. The initial condition for each trajectory is drawn independently from $\mu$. For the $m$-th trajectory with initial condition $\mathbf{x}_m$, there are $T_m$ samples drawn independently from measure $\kappa(\cdot, \mathbf{x}_{m,1})$ where the time elapse of $j$-th sample away from initial condition $\mathbf{x}_{m,1}$ is $t_{m,j}$ ($t_{m,1} =0$), where $m=1,\ldots,M$. Then, we can define the corresponding empirical measure from~\cref{eq:markov_kernel}, 
\begin{equation}
\hat{\nu}_{M, \widehat{T}} = \frac{1}{M} \sum_{m=1}^M \delta_{\mathbf{x}_m} \left( \frac{1}{T_m} \sum_{j=1}^{T_m}\delta_{t_{m,j}} \right),
\end{equation}
where $\displaystyle \widehat{T} = \max_{m=1,\ldots,M} T_m$, and $\hat{\nu}_{M, \widehat{T}}$ uniformly converges to $\nu$ as $M,\widehat{T} \rightarrow \infty$.

Similarly, we can rewrite \cref{eq:recurrent_objective} as the following, 
\begin{align}
\label{eq:finite_recurrent}
&\widehat{\mathbf{W}}_{\mathbf{\Phi}}, \widehat{\mathbf{W}}_{\mathbf{\Psi}}, \widehat{\mathbf{K}} = \lim_{M,T \rightarrow \infty} \widehat{\mathbf{W}}_{\mathbf{\Phi},M,\widehat{T}}, \widehat{\mathbf{W}}_{\mathbf{\Psi},M,\widehat{T}}, \widehat{\mathbf{K}}_{M,\widehat{T}} \\
& \nonumber = \lim_{M,\widehat{T} \rightarrow \infty} \argmin_{\substack{\mathbf{W}_{\mathbf{\Phi}}, \mathbf{W}_{\mathbf{\Psi}}, \mathbf{K}}} \widetilde{\mathcal{J}}_{r,M,\widehat{T}}[\mathbf{\Phi}(\cdot; \mathbf{W}_{\mathbf{\Phi}}), \mathbf{K}] +  \widehat{\mathcal{P}}_{M,\widehat{T}}[\mathbf{\Phi}(\cdot; \mathbf{W}_{\mathbf{\Phi}}), \mathbf{\Psi}(\cdot; \mathbf{W}_{{\mathbf{\Psi}}}), \mathbf{K}] \\
& \nonumber = \lim_{M,\widehat{T} \rightarrow \infty} \argmin_{\substack{\mathbf{W}_{\mathbf{\Phi}}, \mathbf{W}_{\mathbf{\Psi}}, \mathbf{K}}}  \left\lVert \mathbf{\Phi} e^{t\mathbf{K}} - \mathcal{K}_t \mathbf{\Phi} \right\rVert^2_{\widehat{\mathcal{G}}_{M,\widehat{T}}}  + 
\left\lVert
\mathbf{\Psi} \circ \mathbf{\Phi} e^{t \mathbf{K}}  - \mathcal{K}_t \mathcal{I}
\right\rVert^2_{\widehat{\mathcal{G}}_{M, \widehat{T}}^N}
\\
& \nonumber= \lim_{M,\widehat{T} \rightarrow \infty} \argmin_{\substack{\mathbf{W}_{\mathbf{\Phi}}, \mathbf{W}_{\mathbf{\Psi}}, \mathbf{K}}} \frac{1}{M} \sum_{m=1}^{M} \frac{1}{T_m} \bigg(    \sum_{j=2}^{T_m} \lVert \mathbf{\Phi}(\mathbf{x}_{m,1}; \mathbf{W}_{\mathbf{\Phi}}) e^{t_{m,j} \mathbf{K}} -  \mathbf{\Phi} (\mathbf{x}_{m,j}) \rVert^2  \\
& \nonumber \qquad \qquad \qquad\qquad + 
\sum_{j=1}^{T_m}
\left\lVert
\mathbf{\Psi} ( \mathbf{\Phi}(\mathbf{x}_{m,1}; \mathbf{W}_{\mathbf{\Phi}})  e^{t_{m,j} \mathbf{K}}; \mathbf{W}_{\mathbf{\Psi}})  - \mathbf{x}_{m,j} 
\right\rVert^2 \bigg),
\end{align}
where $\widehat{G} = L^2(\mathcal{U}, \widehat{\nu}_{M, \widehat{T}})$.
 Note that \cref{eq:finite_recurrent} generalizes the  LRAN model with a simpler loss function~\cite{otto2017linearly} and the discrete spectrum model in the paper of Lusch et al.~\cite{lusch2017deep}. 

Note that the above data setup also generalizes to cases where data along a single long trajectory  is ``Hankelized"~\cite{otto2017linearly,arbabi2017ergodic}, i.e., dividing a long trajectory into several smaller-sized consecutive trajectories. This would lead to a truncated time horizon in the model, which leads to better computational efficiency compared to the original data, at the cost of some loss in prediction accuracy.

\subsection{Guaranteed stabilization of the Koopman operator}

Eigenvalues of the Koopman operator are critically important in understanding the temporal behavior of certain modes in the dynamical system. For a measure-preserving system~\cite{budivsic2012applied} or systems on an attractor, e.g., post-transient flow dynamics~\cite{arbabi2017study}, even if the system is chaotic, the eigenspectrum of the continuous-time Koopman operator would still be on the imaginary axis. It should be noted that although the  Koopman operator still accepts unstable modes, i.e., the real part of the eigenvalues of continuous-time Koopman operator being positive,  its absence in  systems governed by Navier-Stokes equation in fluid mechanics has been documented~\cite{mezic2013analysis}. Hence, in this work, we assume that the Koopman eigenvalues corresponding to the finite dimensional Koopman invariant subspace of interest have non-positive real parts. It is important to note that the concept of unstable Koopman modes should not be confused with that of flow instability.  Prior models (for instance, ~\cite{otto2017linearly}) have not explicitly taken  stability into account, and thus resulted in slightly unstable Koopman modes. While this is acceptable for relatively short time predictions, long time predictions will be problematic. 

Before presenting the stabilization technique, it is instructive to note the non-uniqueness of \emph{ideal} observation functions $\mathbf{\Phi}$, i.e., the one corresponding to the exact Koopman invariant subspace. This is because one can simultaneously multiply any $D \times D$ invertible real matrix and its inverse before and after the observation vector $\mathbf{\Phi}$ while keeping the Koopman eigenfunction, eigenvalues, and the output from the reconstruction the same. Thus, observation functions described by neural networks cannot be expected to be uniquely determined. We  will leverage this non-uniqueness to enforce stability. 

Enlightened by recent studies~\cite{haber2017stable,chang2018reversible,chang2019antisymmetricrnn} in the design of stable deep neural network structures where skew-symmetric weights are used inside nonlinear activations, we devised a novel parameterization for the realization of the Koopman operator in the following form:
\begin{equation}
\label{eq:stable_structure}
\mathbf{K}_{stable} = 
\begin{bmatrix}
-\sigma^2_1 &       \zeta_1    &            &          \\
  -\zeta_1 &      \ddots    &  \ddots    &          \\
           &        \ddots    &  \ddots   & \zeta_{D-1}  \\
           &                  &   -\zeta_{D-1} &  -\sigma^2_D 
\end{bmatrix},
\end{equation}
where $\zeta_1,\ldots,\zeta_{D-1}$, $\sigma_1,\ldots,\sigma_D \in \mathbb{R}$.

In  \cref{apdx}, we show that the real parts of the eigenvalues of a $n \times n$ real (possibly non-symmetric) negative semi-definite matrix $\mathbf{A}$ are non-positive. We then prove that the above parameterization posits a constraint that the time evolution associated with $\mathbf{K}_{stable}$ in \cref{eq:stable_structure} for any choice of parameters in $\mathbf{R}$ is always stable. Further, the constraint from the parameterization in \cref{eq:stable_structure} is actually rich enough such that any diagonalizable matrix corresponding to a stable Koopman operator can be represented without loss of expressivity.
\begin{theorem}
For any real square diagonalizable matrix $\mathbf{K} \in \mathbb{R}^{D\times D}$ that only has non-positive real parts of the eigenvalues $D \ge 2$, there exists a set of $ \zeta_1,\ldots,\zeta_{D-1}$, $\sigma_1,\ldots,\sigma_D \in \mathbb{R}$ such that the corresponding $\mathbf{K}_{stable}$ in \cref{eq:stable_structure} is similar to $\mathbf{K}$ over $\mathbb{R}$. Moreover, for any $\zeta_1,\ldots,\zeta_{D-1}$, $\sigma_1,\ldots,\sigma_D \in \mathbb{R}$,  the real part of the eigenvalues of the corresponding $\mathbf{K}_{stable}$ is non-positive. 
\end{theorem}
\begin{proof}
For any real square diagonalizable matrix $\mathbf{K} \in \mathbb{R}^{D\times D}$ that only has non-positive real parts of the eigenvalues, there exists an eigendecomposition, 
\begin{equation}
\mathbf{K} = \mathbf{M} \mathbf{J} \mathbf{M}^{-1}, \quad \mathbf{M}, \mathbf{J} \in \mathbb{C}^{D \times D}.
\end{equation}
Without loss of generality, the diagonal matrix $\mathbf{J}$ contains $2D_c$ complex eigenvalues $\{\lambda^c_j\}_{j=1}^{2D_c}$ and $D_r$ real eigenvalues $\{\lambda^r_j\}_{j=1}^{D_r}$ where  $2D_c + D_r = D$. 

Consider a $2\times 2$ real matrix, $\mathbf{A}_j = \begin{bmatrix}
-\sigma_j^2 & \zeta_j\\
- \zeta_j & -\sigma_{j}^2
\end{bmatrix}$ where the eigenvalues are $\lambda_{1,2} = - \sigma_j^2 \pm j\zeta_j^2$. We use this matrix to construct a $2\times 2$ matrix that has   eigenvalues $\lambda^c_j$.  For each pair of complex eigenvalues, $j=1,\ldots,D_c$, we have $\mathbf{A}_j = \begin{bmatrix}
\textrm{Re}(\lambda^c_j) & \sqrt{  |\textrm{Im}(\lambda^c_j)|  } \\
- \sqrt{ |\textrm{Im}(\lambda^c_j)| } & \textrm{Re}(\lambda^c_j)
\end{bmatrix}$. Next, combining with the $D_r$ real eigenvalues, we have the following block diagonal matrix that shares the same eigenvalues as $\mathbf{J}$, 
\begin{equation}
\widetilde{\mathbf{K}}=
\begin{bmatrix}
\mathbf{A}_1    &                &                   &             &            &           \\ 
                & \ddots         &                   &             &            &           \\
                &                &  \mathbf{A}_{D_c} &             &            &           \\
                &                &                   & \lambda_1^r &            &           \\
                &                &                   &             & \ddots     &           \\
                &                &                   &             &            &  \lambda^r_{D_r}     \\            
\end{bmatrix} \in \mathbb{R}^{D \times D}.
\end{equation}
Since  two matrices with the same (complex) eigenvalues are similar~\cite{meyer2000matrix}, $\mathbf{K}$ is similar to $\widetilde{\mathbf{K}}$ over $\mathbb{C}$. However, since  $\mathbf{K}$ and $\widetilde{\mathbf{K}}$ are matrices over $\mathbb{R}$,  they are already similar over $\mathbb{R}$ (Refer Corollary 1 on p. 312 in Ref.~\cite{Herstein}). Notice that the general form of $\widetilde{\mathbf{K}}$ is just a special form of $\mathbf{K}_{stable}$ in \cref{eq:stable_structure}. Therefore, one can always find a set of $\zeta_1,\ldots,\zeta_D$ and $\sigma_1,\ldots,\sigma_D$ such that $\mathbf{K}_{stable}$ is similar to $\mathbf{K}$ over $\mathbb{R}$. Next, notice that given $\sigma_1,\ldots,\sigma_D$, $\zeta_1,\ldots,\zeta_D$, $\forall \mathbf{v}\in \mathbb{R}^{D \times 1}$, we have
\begin{equation}
\mathbf{v}^\top \mathbf{K}_{stable} \mathbf{v} = \frac{1}{2}\mathbf{v}^\top (\mathbf{K}_{stable} +\mathbf{K}_{stable}^\top) \mathbf{v} = \mathbf{v}^\top \begin{bmatrix}
-\sigma_1^2    &      &     \\
               & \ddots & \\
               & & -\sigma_D^2 
\end{bmatrix} \mathbf{v} \le 0.
\end{equation}
Therefore, $\mathbf{K}_{stable}$ is a negative semi-definite matrix. Following \cref{lem:lemma_negative_semi}, the real part of any eigenvalue of $\mathbf{K}_{stable}$ is non-positive.
\end{proof}

Thanks to the conjugation over $\mathbb{R}$, the last layer of encoder, i.e., $\mathbf{\Phi}, $ and the first layer of the decoder, i.e., $\mathbf{\Psi}$, will absorb any linear transformation necessary.
Hence, the above theorem simply means one can parameterize $\mathbf{K}$ with \cref{eq:stable_structure} that guarantees stability of the linear dynamics of $\mathbf{\Phi}$ without loss of expressibility for cases where unstable modes are absent.  

Careful readers might notice that indeed one can further truncate the parameterization with a $2\times 2$ block matrix instead of the tridiagonal form used in \cref{eq:stable_structure}. However, since this reduction takes more effort in the implementation while the parameterization in \cref{eq:stable_structure} has already reduced the number of parameters from $O(D^2)$ to $O(D)$, we prefer the tridiagonal form. For the rest of the work, we will use this parameterization for all of the cases concerned.

\subsection{Design of neural network architecture with SVD-DMD embedding}

In the previous subsection, we described the deep learning formulation for the Koopman operator as a finite dimensional optimization problem that approximates the constrained variational problem, together with a stable parameterization for $\mathbf{K}$. In this subsection, we will describe the design of our neural network that further embeds the SVD-DMD for differential and recurrent forms,  which were presented  in \cref{sec:diff_form,sec:recurrent}.

\subsubsection{Input normalization}
\label{sec:input_norm}

As a standard procedure, we consider normalization on the snapshot matrix of state variable $\mathbf{X}$,
\begin{equation}
\label{eq:data}
\mathbf{X} = \begin{bmatrix}
\mathbf{x}_1 \\ 
\vdots \\
\mathbf{x}_M 
\end{bmatrix} \in \mathbb{R}^{M \times N}
\end{equation}
for better training performance on neural networks~\cite{goodfellow2016deep}. Specifically, we consider Z-normalization shown in \cref{eq:z_n}, i.e., subtracting the mean of $\mathbf{x}$ then dividing the standard deviation to obtain $\mathbf{z}$. 
\begin{equation}
\label{eq:z_n}
\mathbf{z} = (\mathbf{x} - \overline{\mathbf{x}}) \mathbf{\Lambda}^{-1},
\end{equation}
where $\overline{\mathbf{x}} = \frac{1}{M} \sum_{m=1}^{M} \mathbf{x}_m$, $\mathbf{\Lambda} = \textrm{diag}\{ d_1,\ldots,d_N \}$, $d_j$ is the uncorrected standard deviation of $j$-th component of $\mathbf{x}$, i.e., 
\begin{equation}
d_j = \sqrt{\frac{1}{M} \sum_{m=1}^{M} (\mathbf{x}_{m,j} - \overline{\mathbf{x}}_j)^2} ,
\end{equation}
where $\mathbf{x}_{m,j}$ is $j$-th component of $m$-th data, $j=1,\ldots,N$.

While such a normalization is helpful for neural network training in most cases, in some cases where the data is a set of POD coefficients, it cannot differentiate between components that could be more significant than others. Therefore, for those cases, we consider a different normalization with the $\mathbf{\Lambda}$ that sets the ratio of standard deviation between components as: 
\begin{equation}
\label{eq:z_n_a}
\mathbf{\Lambda} = d_{max}\mathbf{I},
\end{equation}
where $d_{max} = \max_{j=1,\ldots,N} d_j$, and $\mathbf{I}$ is the identity matrix.

\subsubsection{Embedding with SVD-DMD}

\label{sec:svd_dmd}


 Instead of directly using the standard feedforward neural network structure employed in previous works, we embed SVD-DMD~\cite{schmid2010dynamic} into the framework and learn the residual. Recall for the standard SVD-DMD algorithm, given $M$ sequential snapshots in \cref{eq:data}  uniformly sampled at intervals  $\Delta t$, one linearly approximates the action of $\mathcal{K}_{\Delta t}$ of the first $r$ dominant SVD modes of \emph{centered} snapshots, i.e., $\mathbf{q} = (\mathbf{x} - \overline{\mathbf{x}} ) \mathbf{V}_r$. Here, $r$ is empirically chosen as a balance between numerical robustness and reconstruction accuracy, and $\mathbf{V}_r$ are the first $r$ columns from $\mathbf{V}$ of the SVD of centered snapshots, 
\begin{equation}
\overline{\mathbf{X}} = \begin{bmatrix}
\mathbf{x}_1 - \overline{\mathbf{x}} \\ 
\vdots \\
\mathbf{x}_M - \overline{\mathbf{x}}
\end{bmatrix} = \mathbf{U} \mathbf{\Sigma} \mathbf{V}^\top.
\end{equation}
Then the SVD-DMD operator is simply the matrix $\mathbf{A}$ that minimizes $\sum_{m=1}^{M-1} \lVert \mathbf{q}_{m+1} - \mathbf{q}_m \mathbf{A} \rVert $, where $\mathbf{q}_m$ is the corresponding orthogonal projection of $\mathbf{x}_m$.

\begin{figure}[hpbt]
  \centering
  \includegraphics[width=0.9\linewidth]{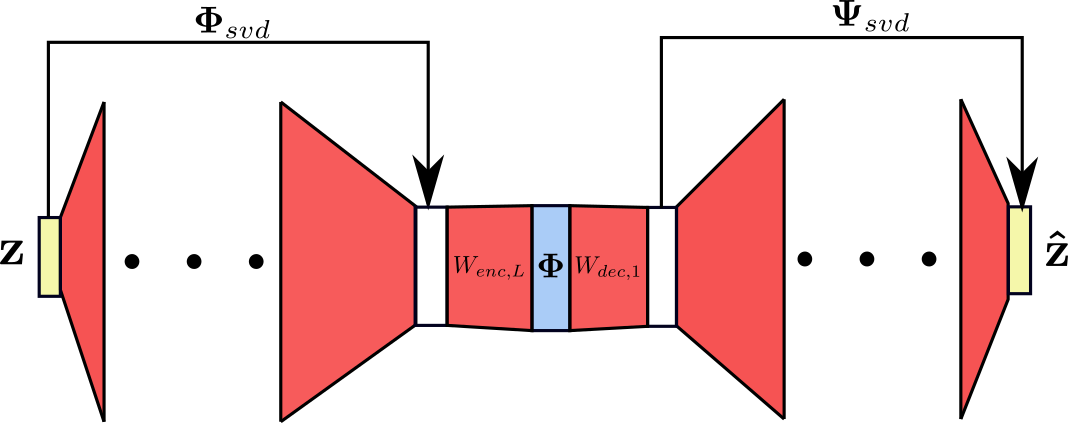}
  \caption{Sketch of SVD-DMD embedding in the feedforward neural network. Red blocks represent neural network weights and biases. Yellow blocks represent states. Light blue block represent observables.}
  \label{fig:svd_dmd_embedding}
\end{figure}

To embed the above SVD-DMD structure into the neural network, we introduce two modifications. First, we take $r=D$. Since $D$ is arbitrary, if $N < D$, we simply append zero columns in $\mathbf{V}_D$. Second, we would also need to accommodate $\mathbf{z}$ with the input normalization in \cref{sec:input_norm}. Thus, we cast $\mathbf{q} = \mathbf{z} \mathbf{\Lambda} \mathbf{V}_D$. We then have,  
\begin{align}
\mathbf{\Phi}_{\textcolor{black}{svd}}(\mathbf{z}) &= \mathbf{z} \mathbf{\Lambda} \mathbf{V}_D, \quad \mathbf{\Psi}_{\textcolor{black}{svd}}(\mathbf{\Phi})\triangleq \mathbf{\Phi}  \mathbf{V}_D^\top \mathbf{\Lambda}^{-1},
\\
\mathbf{\Phi}(\mathbf{z}) &= \textcolor{black}{\underbrace{\mathbf{\Phi}_{nn}(\mathbf{z}) W_{enc,L}  }_{\textrm{nonlinear observables}} + \underbrace{\mathbf{\Phi}_{svd}(\mathbf{z})W_{enc,L}}_{\textrm{linear observables}}  }, \\  
\mathbf{\Psi}(\mathbf{\Phi}) &= \textcolor{black}{\underbrace{\mathbf{\Psi}_{nn}(\mathbf{\Phi})}_{\textrm{nonlinear reconstruction}}  + \underbrace{ \mathbf{\Psi}_{\textcolor{black}{svd}}( \mathbf{\Phi} W_{dec,1})}_{\textrm{linear reconstruction}} }, 
\end{align} 
where $\mathbf{\Phi}_{nn} \triangleq \mathbf{\Phi}(\cdot; \textcolor{black}{\mathbf{W}_\mathbf{\Phi} \setminus{\{W_{enc,L}\}}})$, \textcolor{black}{$\mathbf{W}_\mathbf{\Phi} = \{W_{enc,1}, b_{enc,1},\ldots, W_{enc,L}\}$,} $\mathbf{\Psi}_{nn} \triangleq \mathbf{\Psi}(\cdot; \mathbf{W}_\mathbf{\Psi})$, \textcolor{black}{$\mathbf{W}_{\mathbf{\Psi}} = \{W_{dec,1},\ldots, W_{dec,L}, b_{dec,L} \}$, $L$ is the number of layers for encoder or decoder neural network}. \textcolor{black}{The embedding is illustrated in \cref{fig:svd_dmd_embedding}.}


%

The intuition behind the embedding of the SVD-DMD into the framework is given below:

\begin{enumerate}
\item Schmid's DMD algorithm~\cite{schmid2010dynamic}  relies on the dominant POD modes to reduce the effect of noise, and has been shown to approximate the Koopman invariant subspace~\cite{rowley2009spectral}. This has been demonstrated even for high-dimensional nonlinear dynamical systems with millions of degrees of freedom~\cite{schmid2011applications}. Moreover, for systems with continuous spectra, POD appears to be a robust alternative to Koopman mode decomposition~\cite{mezic2015applications}. Thus, we assume that the true Koopman invariant subspace is easier to obtain by only learning the residual with respect to the dominant POD subspace. 
\item Although we did not precisely implement ResNet blocks~\cite{he2016deep} (empirically, for the problem concerned in this paper, our network does not have to be as deep as common architectures in the deep learning community~\cite{he2016deep}), we believe that such fixed mappings may  have similar benefits in ResNets. 
\item For an ideal linear dynamical system, using the aforementioned neural network model with nonlinear reconstruction can result in an infinite number of unnecessary Koopman modes as global minima. For example, consider the case of a linear dynamical system, i.e., $\mathbf{F}(\mathbf{x}) = \mathbf{x} \mathbf{A} $, $\mathbf{A} \in \mathbb{R}^{N \times N}$. The \emph{desired} Koopman invariant subspace is trivially the span of the projections of $\mathbf{x}$ onto each component and the \emph{desired} Koopman eigenvalues are simply that of $\mathbf{A}$. Assuming simple and real eigenvalues, one can have $\dot{\mathbf{y}} = \mathbf{y}\mathbf{\Lambda}$, where $\mathbf{y} = \mathbf{x} \mathbf{M}$, $\mathbf{A} = \mathbf{M} \mathbf{J} \mathbf{M}^{-1}$ as the eigen-decomposition with $\mathbf{J} = \textrm{diag}\{ \lambda_j \}_{j=1}^{N}$. Then for each component $y_j$, we have $\dot{y_j} = \lambda_j y_j$. For any $n_j \in \mathbb{N}$,  consider the observable $\phi_j = y_j^{2n_j + 1}$, one can have $\dot{\phi_j} = (2n_j + 1) \lambda_j \phi_j$, i.e., $\textrm{span}\{\phi_1,\ldots,\phi_N\}$ is invariant to $\mathcal{K}_t$~\footnote{Budi\v{s}i\'{c} et al.~\cite{budivsic2012applied} showed that the set of eigenfunctions naturally forms an Abelian semigroup under pointwise products.}. Then, consider the nonlinear decoder as one that simply takes the $2n_j+1$-th root on $\phi_j$, and one can recover $\mathbf{y}$ exactly. Finally, augmenting the decoder with $\mathbf{M}^{-1}$ and the encoder with $\mathbf{M}$, the neural network model can find a spurious linear embedding, with eigenvalues as $\{ (2n_1+1)\lambda_1,\ldots, (2n_N+1)\lambda_N\}$ rather than the \emph{desired} $\{\lambda_1,\ldots,\lambda_N\}$, which is an over-complicated nonlinear reconstruction. It is trivial to generalize the above thought experiment to cases where eigenvalues are complex for a real-input-real-output neural network to accommodate. On the other hand,  since most often the neural network is initialized with small weights near zero, the effect of the nonlinear encoder and decoder can be small initially compared to the DMD part. Thus, if the system can be exactly represented by DMD, the optimization for the embedded architecture is initialized near the \emph{desired} minimum. We note that \textcolor{black}{this could lead to  attenuation of} the spurious modes due to the nonlinear reconstruction for essentially linear dynamics.\footnote{It has to be mentioned that that such an issue could exist also in cases where DMD is not desired.}

\end{enumerate}



\begin{figure}[hpbt]
  \centering
  \includegraphics[width=0.5\linewidth]{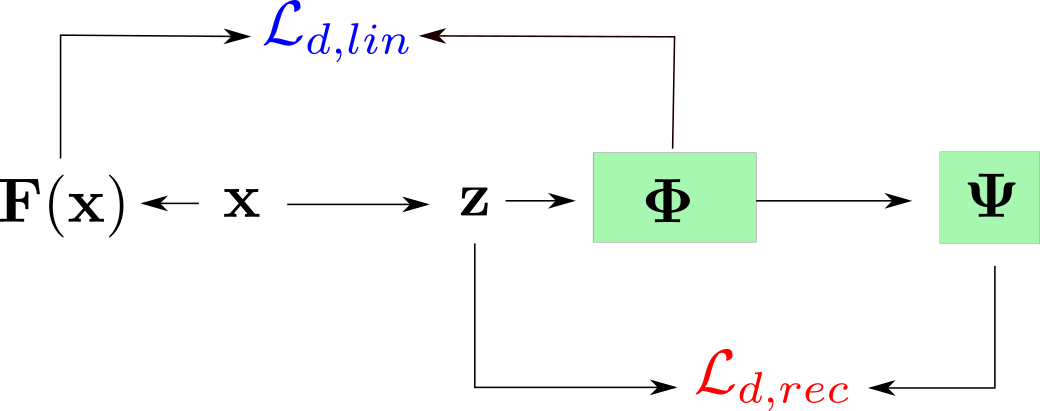}
  \caption{Sketch of the framework of learning a continuous-time Koopman operator in the differential form.}
  \label{fig:framework_diff}
\end{figure}

\begin{figure}[hpbt]
  \centering
  \includegraphics[width=0.8\linewidth]{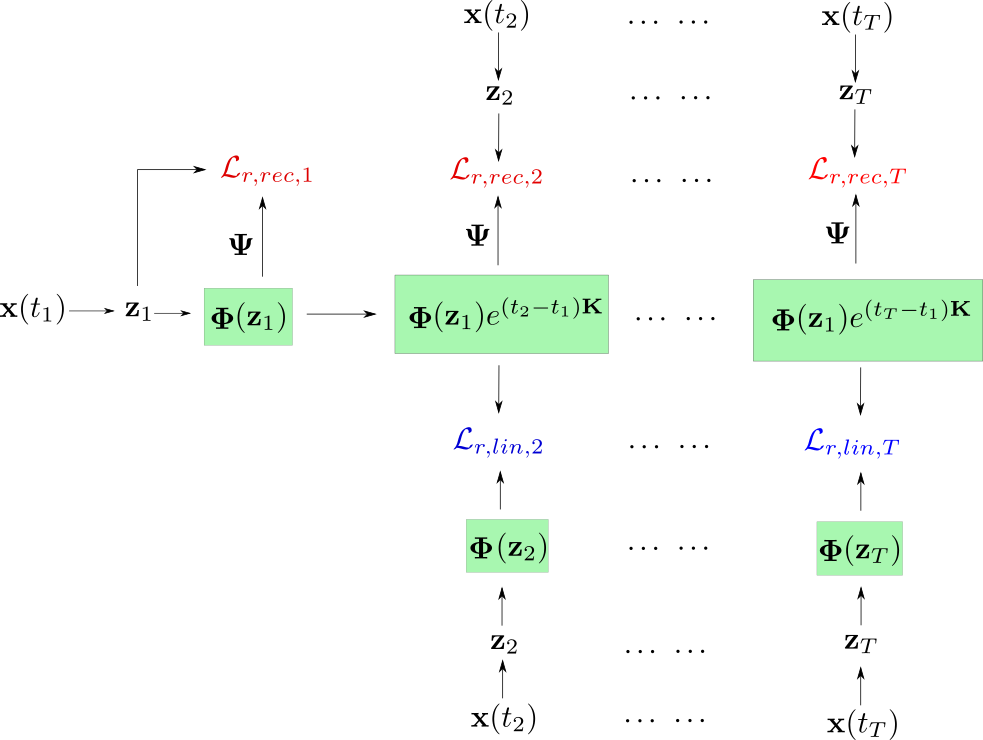}
  \caption{Sketch of the framework of learning a continuous-time Koopman operator in the recurrent form.}
  \label{fig:framework_recurrent}
\end{figure}

The neural network architecture for the differential form  is shown  in \cref{fig:framework_diff},
and for the recurrent form in \cref{fig:framework_recurrent}. 
Note that, if \textcolor{black}{the nonlinear part, i.e., the feedforward} neural network is not activated, the above formulation reverts  to \textcolor{black}{an over-parameterized} SVD-DMD. Specifically, the recurrent form model with neural network deactivated can be viewed as a simplified version of optimized DMD~\cite{askham2018variable,chen2012variants}. 

\subsection{Implementation}
The framework is built using Tensorflow~\cite{tensorflow}. Neural network parameters $\mathbf{W}_{\mathbf{\Phi}}, \mathbf{W}_{\mathbf{\Psi}}$, are initialized with the standard truncated normal distribution.  $\mathbf{K}$ is initialized with the corresponding DMD approximation. The objective function is optimized using Adam~\cite{kingma2014adam}, which is an adaptive first order stochastic optimization method using gradient updates scaled by square roots of exponential moving averages of previous squared gradients.   
Note that we also include weight decay regularization, $\mathcal{L}_{reg} = \left\lVert \bm{W}_{\mathbf{\Phi}} \right\rVert^2 +  \left\lVert \bm{W}_{\mathbf{\Psi}} \right\rVert^2$, in the objective function, to avoid spurious oscillations in the learned Koopman functions, which helps generalization in an interpolation sense~\cite{goodfellow2016deep}.

\section{Probabilistic formulation}

\label{sec:bayesdl}

\subsection{Bayesian neural networks}

A Bayesian formalism is adopted to quantify the impact of several sources of uncertainty in model construction on the model predictions. 
Bayes' rule is
\begin{equation}
{P}( \mathbf{\Theta} | \mathcal{D}) = \frac{{P}(\mathcal{D}|\mathbf{\Theta})   {P}(\mathbf{\Theta}) }{{P}(\mathcal{D})} \iff \textrm{Posterior} = \frac{\textrm{Likelihood}\times \textrm{Prior}}{\textrm{Evidence}},
\end{equation}
where $\mathcal{D}$ is data, and $\mathbf{\Theta}$ is the set of parameters. For simplicity,  $P$ represents the probability density function (PDF) on the measure space generated by the data and parameters. 
 From a traditional Bayesian standpoint, as the number of parameters in the neural network is large, it is impossible for common inference tools such as Markov Chain Monte Carlo (MCMC) to be practical. 
To overcome the curse of dimensionality, several general approaches such as Laplacian approximation~\cite{denker1991transforming} and variational inference~\cite{kucukelbir2017automatic,blundell2015weight}   have been proposed. The former is computationally economical but has two major limitations: First, computing the full Hessian is impossible and expensive for a high dimensional problem and most often   approximations such as $\mathcal{J}^\top \mathcal{J}$, where  $\mathcal{J}$ is the Jacobian are employed \cite{mackay1992practical}. Second, it only provides a local approximation, which can often be far removed from the true posterior. Variational inference has become popular in the deep learning community  as it offers an informed balance between the computationally expensive MCMC method, and  the cheap but less descriptive models such as the Laplacian approximation.  Historically, variational inference for neural networks~\cite{hinton1993keeping} has been difficult \cite{neal2012bayesian} to formulate, largely due to the difficulty of deriving analytical solutions to the required integrals over the variational posteriors \cite{graves2011practical} even for simple network structures. Graves \cite{graves2011practical} proposed a stochastic method for variational inference with a diagonal Gaussian posterior that can be applied to almost any differentiable log-loss parametric model, including neural networks. However,  there is always a trade-off between the complexity of the posterior and scalability and robustness~\cite{vireview}. In this work, we adopt the  mean-field variational inference~\cite{kucukelbir2017automatic}. 



\subsection{Mean-field variational inference}

\begin{figure}[hpbt]
  \centering
  \includegraphics[width=0.8\linewidth]{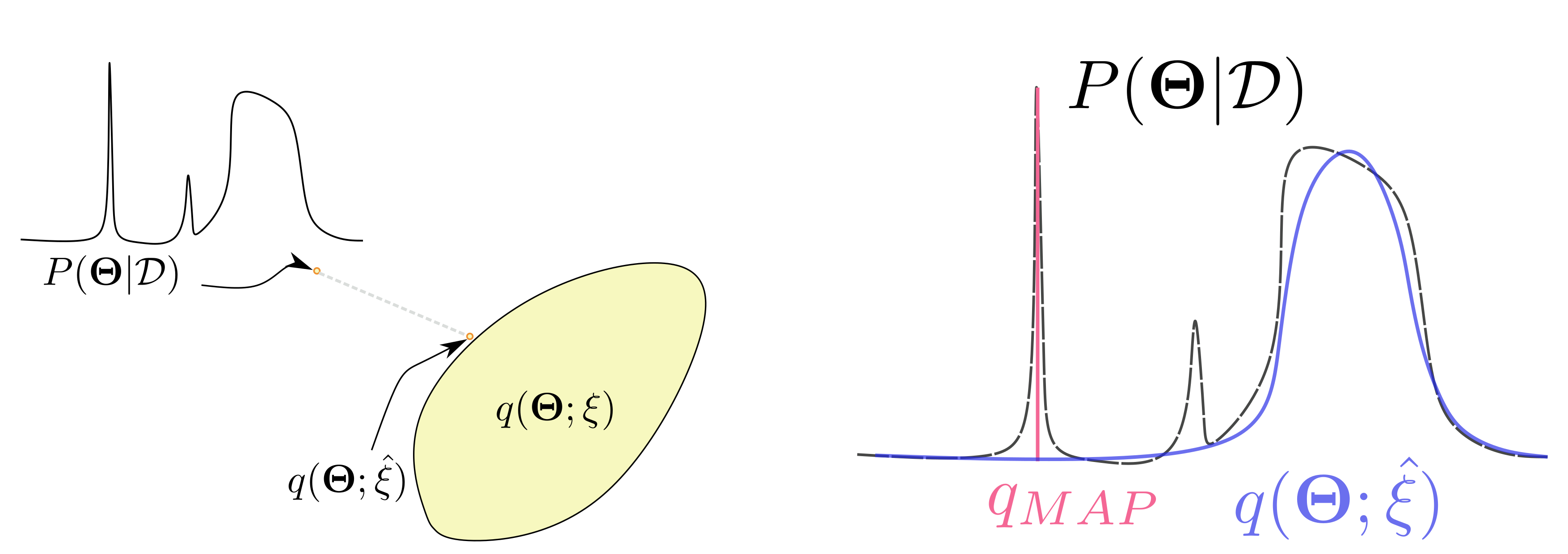}
  \caption{Left: illustration of variational inference. Right: difference between MAP and variational inference.}
  \label{fig:vi}
\end{figure}

As illustrated in the left figure of \cref{fig:vi}, the key idea in variational inference~\cite{blei2017variational} is to recast Bayesian inference as an optimization problem by searching the best parameterized probability density function $q(\mathbf{\Theta};\hat{\xi})$ in a family of approximating densities, namely the variational posterior, $\{ q(\mathbf{\Theta};\xi) \vert \xi \in \Xi \}$, such that it is closest to the true posterior $P(\mathbf{\Theta}|\mathcal{D})$. Most often, the Kullback--Leibler (KL) divergence is employed to measure the distance, which is defined as $\textrm{KL}( q(\mathbf{\Theta}; {\xi}) \Vert P(\mathbf{\Theta}|\mathcal{D}) ) = \int_{\Omega} q(\mathbf{\Theta}; \xi) \log\frac{q(\mathbf{\Theta}; \xi)}{P(\mathbf{\Theta} \vert \mathcal{D})}  d\mathbf{\Theta} = \mathbb{E}_{q(\mathbf{\Theta}; \xi)} \big[ \log{q(\mathbf{\Theta}; \xi)}{P(\mathbf{\Theta} \vert \mathcal{D})} \big] $, where $\Omega$ is the support of $q(\mathbf{\Theta}; \xi)$. This implies $\textrm{supp}(q(\mathbf{\Theta}; \xi)) \subseteq \textrm{supp}(P(\mathbf{\Theta} \vert \mathcal{D}))$. Further, we assume $\textrm{supp}(P(\mathbf{\Theta} \vert \mathcal{D})) = \textrm{supp}(P(\mathbf{\Theta}))$. $\Xi$ is the domain of $\xi$, depending on the parameterization and family of approximating densities. 

However, direct computation of the KL divergence is intractable, since we do not have access to $\log P( \mathcal{D})$. Instead, we choose an alternative, the evidence lower bound (ELBO), i.e., the negative KL divergence plus $\log P(\mathcal{D})$, in \cref{eq:elbo} to be maximized. 
\begin{equation}
\label{eq:elbo}
\mathcal{L}_{elbo} (\xi) = \mathbb{E}_{q(\mathbf{\Theta}; \xi)} \big[ \log P(\mathcal{D}, \mathbf{\Theta}) \big] - \mathbb{E}_{q(\mathbf{\Theta}; \xi)} \big[ q(\mathbf{\Theta}; \xi) \big].
\end{equation}

To maximize the ELBO, we leverage  automatic differentiation from Tensorflow to compute the gradients with respect to $\xi$, following the framework of Automatic Differentiation Variational Inference (ADVI)~\cite{kucukelbir2017automatic} where Gaussian distributions are considered as the variational family. Specifically, we employ the mean-field assumption in \cref{eq:mean_field} such that,
\begin{equation}
\label{eq:mean_field}
q(\mathbf{\Theta}; \xi) = \prod_{j=1}^{Z} q(\mathbf{\theta}_j; \xi_j),
\end{equation}
where $Z$ is the total number of parameters. For $j=1,\ldots,Z$, $\theta_j$ is the $j$-th parameter as random variable, and $\xi_j$ is the corresponding variational parameter that describes the distribution. This is particularly convenient for neural network models constructed in Tensorflow since weights and biases are naturally defined on some real coordinate space. If the support of the parameter distribution is restricted, one can simply consider a one-to-one differentiable coordinate transformation $\Upsilon(\mathbf{\Theta}) = \mathbf{Z}$, such that $\mathbf{Z}$ is not restricted in some real coordinate space, and posit a Gaussian distribution on $\mathbf{Z}$. Note that this naturally induces non-Gaussian distribution. Here, we employ the ADVI functionality in Edward~\cite{tran2016edward}, which is built upon Tensorflow to implement ADVI. Interested readers should refer to the original paper of ADVI~\cite{kucukelbir2017automatic} for the specific details of implementing mean-field variational inference including the usage of the reparameterization-trick to compute the gradients. In contrast, note that the maximum a posteriori (MAP) estimation of the posterior can be cast as a regularized deterministic model as  illustrated in \cref{fig:vi}. Since weight decay is employed in previous deep learning models~\cite{otto2017linearly,lusch2017deep} to learn the Koopman decomposition, one can show that the previous model is essentially a MAP estimation of the corresponding posterior. 

Note that the mean-field Gaussian assumption is still simplified, yet effective and scalable for deep neural nets~\cite{blundell2015weight}. It is interesting to note that several recent works~\cite{zhu2018bayesian,zhu2019physics} leverage Stein Variational Gradient Descent (SVGD)~\cite{NIPS2016_6338}, a non-parametric, particle-based inference method, which is able to capture multi-modal posteriors. Robustness and scalability to high dimensions, e.g., deep neural nets, is still an area of active research~\cite{2017arXiv171104425Z,liu2018stein,shi2018kernel,lu2019scaling}.

\subsection{Bayesian hierarchical model setup}

Recall that we have the following parameters for the deep learning models introduced in \cref{sec:diff_form} and \cref{sec:recurrent}: 1) weights and biases for the ``encoder", $\mathbf{W}_{\mathbf{\Phi}}$, 2) weights and biases for the ``decoder", $\mathbf{W}_{\mathbf{\Psi}}$, and 3) stabilized $\mathbf{K}$ with $\zeta_1,\ldots,\zeta_D$ and $\sigma_1,\ldots,\sigma_D$. Based on mean-field assumptions, we just need to prescribe the prior and variational posterior for each parameter.
For each weight and bias, we posit a Gaussian prior with zero mean, with the scaling  associated with each parameter to follow the recommended half--Cauchy distribution~\cite{gelman2004prior,polson2012on}, which has zero mean and scale as 1 (empirical). The variational posterior for each weight and bias is Gaussian, and Log-normal  for scale parameters. For the off-diagonal part of $\mathbf{K}$, $\zeta_i$, we posit the same type of Gaussian prior as before with the scale parameter following a hierarchical model for each $i=1,\ldots, D$. The corresponding variational posterior for $\zeta_i$ is still Gaussian while log-normal for the scale parameter. For the non-negative diagonal part of $\mathbf{K}$, we posit a Gamma distribution for each $\sigma_i^2$, $i=1,\ldots,D$, with rate parameter as 0.5 and shape parameter following the previous half-Cauchy distribution. Variational posteriors for both $\sigma^2_i$ and its shape parameter are log-normal. 
%
%

%
%
For the differential form in \cref{sec:diff_form}, we set up the following likelihood function\footnote{Note that independence between the data and the structure of the aleatoric uncertainty is assumed. Such an assumption correlates well with existing deterministic models based on mean-square-error. } based on Z-normalized data $\mathcal{D} = \{\mathbf{z}_{m}, \mathbf{\dot{z}}_m\}_{m=1}^{M}$:
\begin{equation}
\label{eq:likeli_diff}
\mathcal{D} | \mathbf{\Theta} \sim \prod_{m=1}^{M}
\mathcal{N}\left(
\begin{bmatrix}
\mathbf{\Psi}(\mathbf{\Phi}(\mathbf{z}_m; \mathbf{W}_{\mathbf{\Phi}}); \mathbf{W}_{\mathbf{\Psi}})-\mathbf{z}_m  \\
\mathbf{\dot{z}}_m \cdot \nabla_{\mathbf{z}} \mathbf{\Phi}(\mathbf{z}_m; \mathbf{W}_{\mathbf{\Phi}}) - \mathbf{\Phi}(\mathbf{z}_m; \mathbf{W}_{\mathbf{\Phi}}) \mathbf{K}_{stable}
\end{bmatrix}
;
\begin{bmatrix}
\mathbf{0} \\
\mathbf{0}
\end{bmatrix},
\begin{bmatrix}
\mathbf{\Lambda}_{rec} &  \\
 & \mathbf{\Lambda}_{lin}
\end{bmatrix}\right),
\end{equation}
where $\mathbf{\Lambda}_{rec}, \mathbf{\Lambda}_{lin}$ are diagonal covariance matrices with the prior of each diagonal element following the previous half--Cauchy distribution. We also posit log-normals to infer the posterior of $\mathbf{\Lambda}_{rec}, \mathbf{\Lambda}_{lin}$. We denote $\mathbf{\Theta}$ as  $\mathbf{W}_{\mathbf{\Phi}},\mathbf{W}_{\mathbf{\Psi}}$, and associated scale parameters together with $\mathbf{\Lambda}_{rec}$, $\mathbf{\Lambda}_{lin}$. 

For the recurrent form in \cref{sec:recurrent}, given normalized data, \begin{equation}
    \mathcal{D} = \{ \{ \{\mathbf{z}_{m,j}\}_{j=1}^{T_m} \}_{m=1}^{M}, \{ \{t_{m,j}\}_{j=1}^{T_m} \}_{m=1}^{M} \},
\end{equation}
we consider the following likelihood, \footnote{Alternative likelihoods can be chosen to account for the variation of aleatoric noise in time, which is well-suited for short-horizon forecasting. However, we are more interested in a free-run situation where only the initial condition is given.}:
\begin{equation}
\label{eq:likeli_recurrent}
\mathcal{D} | \mathbf{\Theta} 
\sim 
\prod_{m=1}^{M}
\prod_{j=1}^{T_m}
\mathcal{N}\left(
\begin{bmatrix}
\mathbf{\Psi}(\mathbf{\Phi} (\mathbf{z}_{m,1}; \mathbf{W}_{\mathbf{\Phi}}) e^{t_{m,j} \mathbf{K}_{stable}}; \mathbf{W}_{\mathbf{\Psi}}) - \mathbf{z}_{m,j}  \\
\mathbf{\Phi}(\mathbf{z}_{m,1}; \mathbf{W}_{\mathbf{\Phi}}) e^{t_{m,j} \mathbf{K}_{stable}} - \mathbf{\Phi}(\mathbf{z}_{m,j}; \mathbf{W}_{\mathbf{\Phi}})
\end{bmatrix}
;
\begin{bmatrix}
\mathbf{0} \\
\mathbf{0}
\end{bmatrix},
\begin{bmatrix}
\mathbf{\Lambda}_{rec} &  \\
 & \mathbf{\Lambda}_{lin}
\end{bmatrix}\right).
\end{equation}


%

%
%
%
%
%
%

\subsection{Propagation of uncertainties}

\label{sec:monte_carlo}

Given data $\mathcal{D}$ and the inferred posterior $P(\mathbf{\Theta} \vert \mathcal{D})$, assuming a noise-free initial condition $\mathbf{z}_0$, we are interested in future state predictions with uncertainties. For the differential form in \cref{eq:likeli_diff}, we have
\begin{align}
\label{eq:prob_diff}
& P(\mathbf{z}(t) \vert \mathbf{z}_0, \mathcal{D}) = \int P(\mathbf{z}_t \vert \mathbf{W}_{\mathbf{\Phi}}, \mathcal{D}, \mathbf{z}_0)  P(\mathbf{W}_{\mathbf{\Phi}} \vert  \mathcal{D} )) d\mathbf{W}_{\mathbf{\Phi}}, \\
\nonumber & = \int P(\mathbf{z}_t \vert \mathbf{\Phi}(\mathbf{z}_0; \mathbf{W}_{\mathbf{\Phi}}), \mathbf{W}_{\mathbf{\Phi}}, \mathcal{D}, \mathbf{z}_0)   P(\mathbf{W}_{\mathbf{\Phi}} \vert  \mathcal{D} )) d\mathbf{W}_{\mathbf{\Phi}}, \\
\nonumber &= \iint  P(\mathbf{z}_t \vert \mathbf{\Phi}(t), \mathbf{W}_{\mathbf{\Phi}}, \mathcal{D}, \mathbf{z}_0 ) P( \mathbf{\Phi}(t) \vert  \mathbf{\Phi}(\mathbf{z}_0; \mathbf{W}_{\mathbf{\Phi}}), \mathbf{W}_{\mathbf{\Phi}}, \mathcal{D}, \mathbf{z}_0)   P(\mathbf{W}_{\mathbf{\Phi}} \vert  \mathcal{D} )) d \mathbf{\Phi}(t)  d\mathbf{W}_{\mathbf{\Phi}} , \\
\nonumber &= \iiint  
P(\mathbf{z}_t \vert \mathbf{W}_{\mathbf{\Psi}}, \mathbf{\Phi}(t) ) 
P( \mathbf{\Phi}(t) \vert \mathbf{\Phi}(\mathbf{z}_0; \mathbf{W}_{\mathbf{\Phi}}))   
P(\mathbf{W}_{\mathbf{\Phi}} \vert  \mathcal{D} )) P(\mathbf{W}_{\mathbf{\Psi}} \vert \mathcal{D})  d \mathbf{\Phi}(t) d \mathbf{W}_{\mathbf{\Psi}} d\mathbf{W}_{\mathbf{\Phi}}, \\
\nonumber &= \iiiint\!\!\!\iint  
P(\mathbf{z}_t \vert \mathbf{\Lambda}_{rec},\mathbf{W}_{\mathbf{\Psi}}, \mathbf{\Phi}(t) ) 
P( \mathbf{\Phi}(t) \vert \mathbf{\Lambda}_{lin}, \mathbf{K}, \mathbf{\Phi}(\mathbf{z}_0; \mathbf{W}_{\mathbf{\Phi}}))   
P(\mathbf{W}_{\mathbf{\Phi}} \vert  \mathcal{D} ))  P(\mathbf{K} \vert \mathcal{D})   \\ 
& \nonumber \qquad\qquad  \qquad  P(\mathbf{W}_{\mathbf{\Psi}} \vert \mathcal{D}) P(\mathbf{\Lambda}_{lin} \vert \mathcal{D})
P(\mathbf{\Lambda}_{rec} \vert \mathcal{D})
d\mathbf{\Lambda}_{lin} d\mathbf{\Lambda}_{rec} d \mathbf{\Phi}(t)  d \mathbf{W}_{\mathbf{\Psi}} d\mathbf{W}_{\mathbf{\Phi}} d\mathbf{K}.
\end{align}

However, $P( \mathbf{\Phi}(t) \vert \mathbf{K}, \mathbf{\Phi}(\mathbf{z}_0; \mathbf{W}_{\mathbf{\Phi}}))$ is unknown because the differential form does not use trajectory information. If we assume multivariate Gaussian white noise with the same covariance in the linear loss $\mathbf{\Lambda}_{lin}$, then one can forward propagate the  aleatoric uncertainty associated with the likelihood function of the linear loss.  Then one immediately recognizes that the continuous-time random process of $\mathbf{\Phi}(t)$ becomes a multivariate Ornstein--Uhlenbeck process, 
\begin{equation}
\label{eq:dynamics_phi_bayesian}
d\mathbf{\Phi}^\top(t) = \mathbf{K}^\top \mathbf{\Phi}^\top(t) dt + \mathbf{\Lambda}_{lin}^{1/2} d \mathbf{B}(t),
\end{equation}
where $\mathbf{B}(t)$ is a $D$-dimensional Gaussian white noise vector with unit variance. Note that~\cite{ross1996stochastic}
\begin{equation}
\label{eq:prob_phit}
\mathbf{\Phi}(t) \vert \mathbf{z}_0, \mathbf{\Lambda}_{lin}, \mathbf{W}_{\mathbf{\Phi}}, \mathbf{K} 
\sim 
\mathcal{N}( \mathbf{\Phi}(\mathbf{z}_0; \mathbf{W}_{\mathbf{\Phi}}) e^{t\mathbf{K}} , \int_0^{t} e^{s\mathbf{K}} \mathbf{\Lambda}_{lin}  e^{s\mathbf{K}^\top} ds  ),
\end{equation}
where $\int_0^{t} e^{s\mathbf{K}} \mathbf{\Lambda}_{lin}  e^{s\mathbf{K}^\top} ds  = \mathbf{vec}^{-1}( -(\mathbf{K} \oplus \mathbf{K})^{-1} (\mathbf{I} - e^{t(\mathbf{K} \oplus \mathbf{K}) }) \mathbf{vec}(\mathbf{\Lambda}_{lin}) )$. $\mathbf{vec}(\cdot)$ is the stack operator, and $\oplus$ is the Kronecker sum~\cite{meucci2009review}. 

It is interesting to note that, since $\mathbf{K}$ is restricted by \cref{eq:stable_structure} and does not contain any eigenvalues with positive real part, the variance in \cref{eq:prob_phit} will not diverge in finite time. One can thus simply draw samples of $\mathbf{\Phi}(t)$ from \cref{eq:prob_phit}. Thus, we approximate \cref{eq:prob_diff} with Monte Carlo sampling from the corresponding variational posterior, 
\begin{equation}
P(\mathbf{z}(t) \vert \mathbf{z}_0, \mathcal{D}) \approx \frac{1}{N_{mc} M_{mc}}
\sum_{i=1}^{N_{mc}}  \sum_{j=1}^{M_{mc}} P(\mathbf{z}_t \vert \mathbf{\Lambda}_{lin}^{(i)}, \mathbf{\Lambda}_{rec}^{(i)}, \mathbf{W}^{(i)}_{\mathbf{\Phi}}, \mathbf{K}^{(i)}, \mathbf{W}^{(i)}_{\mathbf{\Psi}}, \mathbf{\Phi}^{(j)}(t), \mathbf{z}_0 ),  
\end{equation}
where the superscript with parentheses represents the index of samples, $N_{mc}$, $M_{mc}$ are the number of samples corresponding to  variational posteriors and the Ornstein--Uhlenbeck process.

For the recurrent form in \cref{eq:likeli_recurrent}, the posterior predictive distribution of $\mathbf{z}(t)$ given the initial condition is straightforward:
\begin{align}
\label{eq:prob_recurrent}
 P(\mathbf{z}(t) \vert \mathbf{z}_0, \mathcal{D}) = 
\iiiint & P(\mathbf{z}(t) \vert {\mathbf{\Lambda}_{rec}, \mathbf{K}}, \mathbf{W}_{\mathbf{\Psi}}, \mathbf{W}_{\mathbf{\Phi}}, \mathbf{z}_0)  
P(\mathbf{W}_{\mathbf{\Phi}} \vert  \mathcal{D} )
P(\mathbf{W}_{\mathbf{\Psi}} \vert  \mathcal{D} )
P(\mathbf{K} \vert  \mathcal{D} ) \\
& \nonumber
 P(\mathbf{\Lambda}_{rec} \vert \mathcal{D})
d \mathbf{W}_{\mathbf{\Psi}} 
d \mathbf{\Lambda}_{rec}
d \mathbf{K}
d \mathbf{W}_{\mathbf{\Phi}}, \\ \nonumber
\approx  \frac{1}{N_{mc}} &\sum_{i=1}^{N_{mc}} 
 P(\mathbf{z}(t) \vert {\mathbf{\Lambda}^{(i)}_{rec}, \mathbf{K}^{(i)}}, \mathbf{W}^{(i)}_{\mathbf{\Psi}}, \mathbf{W}^{(i)}_{\mathbf{\Phi}}, \mathbf{z}_0).
\end{align}

%
%

\section{Numerical Results \& discussion}
\label{sec:results}
To demonstrate and analyze the approaches presented herein, three numerical examples are pursued.  

\subsection{2D fixed point attractor}
\label{sec:2d_fixed_point_attractor}

Consider the two-dimensional nonlinear dynamical system~\cite{lusch2017deep} with a fixed point,
\begin{align*}
\dot{x}_1 &= \mu x_1,\\
\dot{x}_2 &= \lambda (x_2 - x_1^2),
\end{align*}
where $\mu=-0.05$, $\lambda=-1$.  For this low dimensional system with known governing equations, we consider the differential form in~\cref{sec:diff_form}, with 1600 states as  training samples
, sampled from $[-0.5, 0.5]$ using the standard Latin-Hypercube-Sampling method~\cite{mckay1979comparison} for both $x_1$ and $x_2$. The embedding of SVD-DMD is not employed in this case since it is not very meaningful. The hyperparameters for training are given in  ~\cref{tab:2d_lusch_hyp}.

\begin{table}[htbp]
{   
    \footnotesize
    \caption{Hyperparameters of differential form model for 2D fixed point attractor.}
    \label{tab:2d_lusch_hyp}
    \begin{center}
    \begin{tabular}{|c|c|c|c|c|} \hline
        \bf layer structure & \bf optimizer  & \bf learning rate & \bf total epoch &  \bf batch size \\ \hline
          2-8-16-16-8-2-8-16-16-8-2 &  Adam & 1e-4 & 20000 & 128 \\ \hline
    \end{tabular}
    \end{center}
}
\end{table}



\begin{figure}[htbp]
\centering
\includegraphics[width=0.8\linewidth]{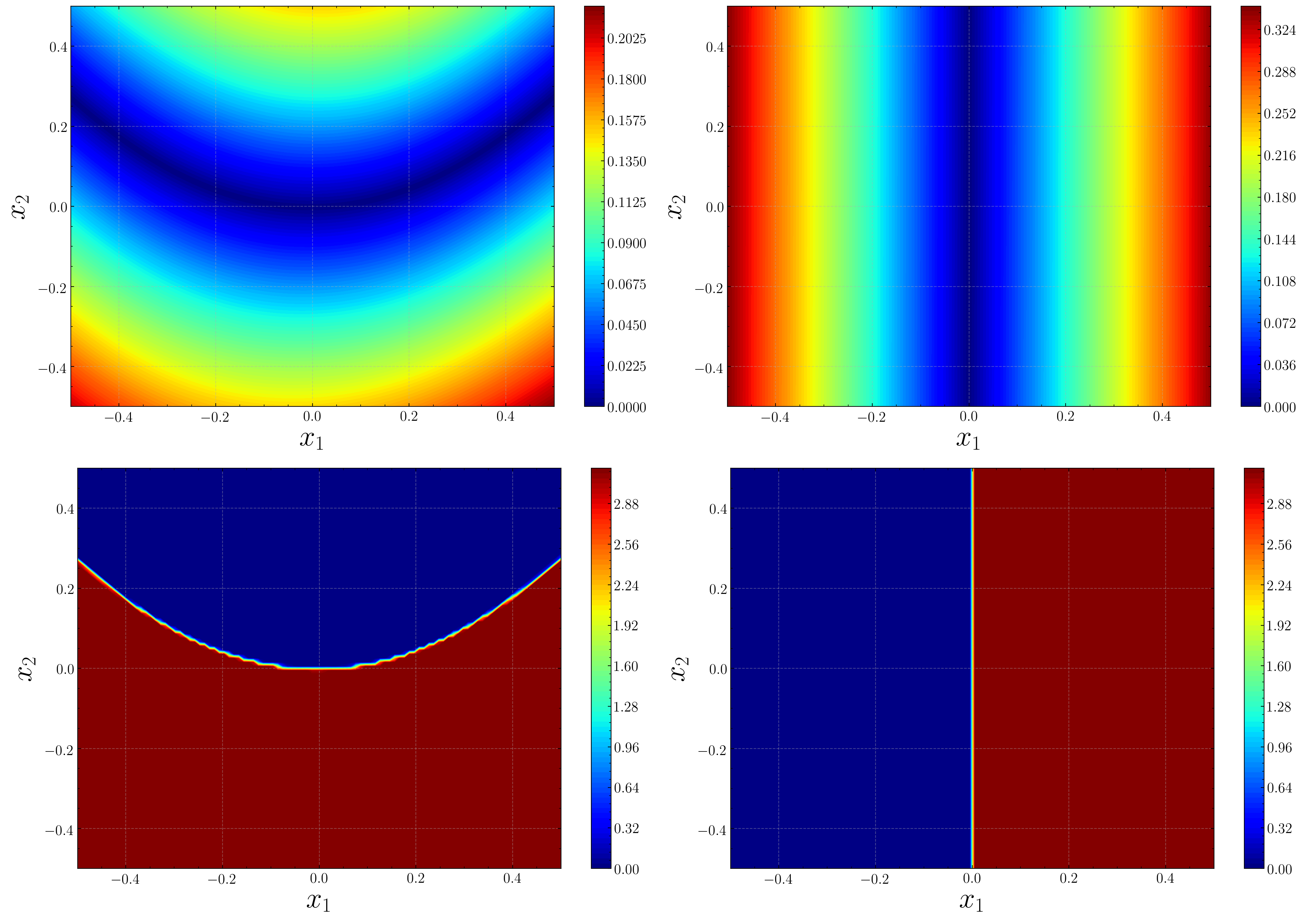}
\label{fig:kf_2d_lusch}
\caption{2D fixed point attractor: mean of the variational posterior of learned Koopman eigenfunctions in differential form. Top left:  mean of the amplitude of Koopman eigenfunction corresponding to $\lambda_1$. Top right:  mean of the amplitude of Koopman eigenfunction corresponding to $\lambda_2$. Bottom left:  mean of the phase angle of Koopman eigenfunction corresponding to $\lambda_1$. Bottom right:  mean of the phase angle of Koopman eigenfunction corresponding to $\lambda_2$. }
\end{figure}

After we obtained the inferred posterior, we consider Monte Carlo sampling described in~\cref{sec:monte_carlo} with $N_{mc}=100$, $M_{mc}=10$ to approximate the posterior distribution of the Koopman observables and prediction on the dynamical system given an initial condition $\mathbf{x}_0$. 
The mean Koopman eigenvalues are $\lambda_1 = -0.99656$ and $\lambda_2 = -0.05049$, and the amplitude and phase angle for the mean eigenfunctions are shown in \cref{fig:kf_2d_lusch}, which resembles the analytic solution with $\varphi_1 =  x_2 - \lambda x_1^2/(\lambda-2\mu)$, $\varphi_2 = x_1$.

\begin{figure}[htbp]
\centering
\includegraphics[width=1\linewidth]{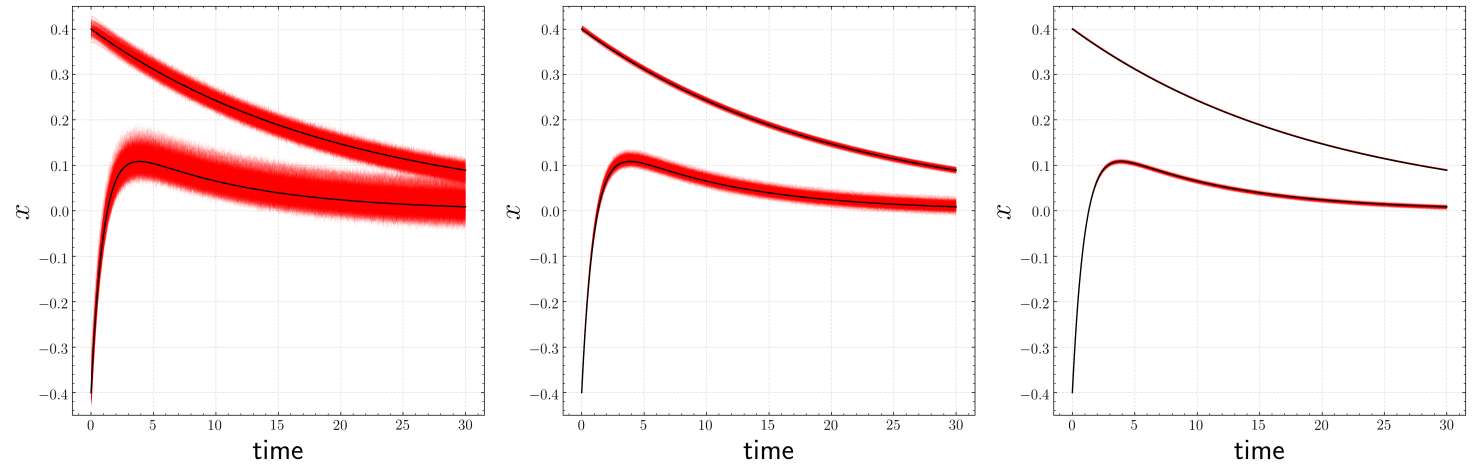}
\label{fig:dynamics_2d_lusch}
\caption{2D fixed point attractor case: Monte Carlo sampling of the predicted trajectory with $\mathbf{x}_0 = [0.4,-0.4]$. Left: 800 (50\% of) original data points. Middle: original 1600 data points. Right: 10000 data points.}
\end{figure}

Finally, since we have a learned a distribution over the Koopman operator, we can obtain the posterior distribution of the predicted   dynamics,  given an arbitrary unseen initial condition following \cref{sec:monte_carlo}, for example, at $\mathbf{x}_0 = (0.4, -0.4)$. The effect of the number of data samples $M$ on the confidence of the predicted dynamics can be seen in \cref{fig:dynamics_2d_lusch}. Clearly, as more data is collected in the region of interest, the propagated uncertainty of the evolution of the dynamics predicted on the testing data decreases as expected. It is interesting to note that, even when the data is halved, the standard deviation of the Koopman eigenvalue is very small compared to the mean. 


\subsection{2D unforced duffing oscillator}

Next, we consider the unforced duffing system:
\begin{align*}
\label{eq:duffing_eq}
\dot{x_1} &=  x_2, \\
\dot{x_2} &= -\delta x_2 - x_1(\beta + \alpha x_1^2),
\end{align*}
where $\delta = 0.5$, $\beta = -1$, $\alpha = 1$. 
We use 10000 samples of the state $\mathbf{x}$ distributed on $x_1,x_2 \in [-2,2]$ from Latin-Hypercube-Sampling. 
We infer the posterior using the differential form model in \cref{sec:diff_form}, with the following hyperparameters in \cref{tab:2d_duffing_hyp}.
\begin{table}[htbp]
{   
    \footnotesize
    \caption{Hyperparameters of differential form model for unforced Duffing system.}
    \label{tab:2d_duffing_hyp}
    \begin{center}
    \begin{tabular}{|c|c|c|c|c|} \hline
        \bf layer structure & \bf optimizer  & \bf learning rate & \bf total epoch &  \bf batch size \\ \hline
          2-16-16-24-16-16-3-16-16-24-16-16-2 &  Adam & 1e-3 & 20000 & 128 \\ \hline
    \end{tabular}
    \end{center}
}
\end{table}

\begin{figure}[htbp]
\centering
    \includegraphics[width=1\linewidth]{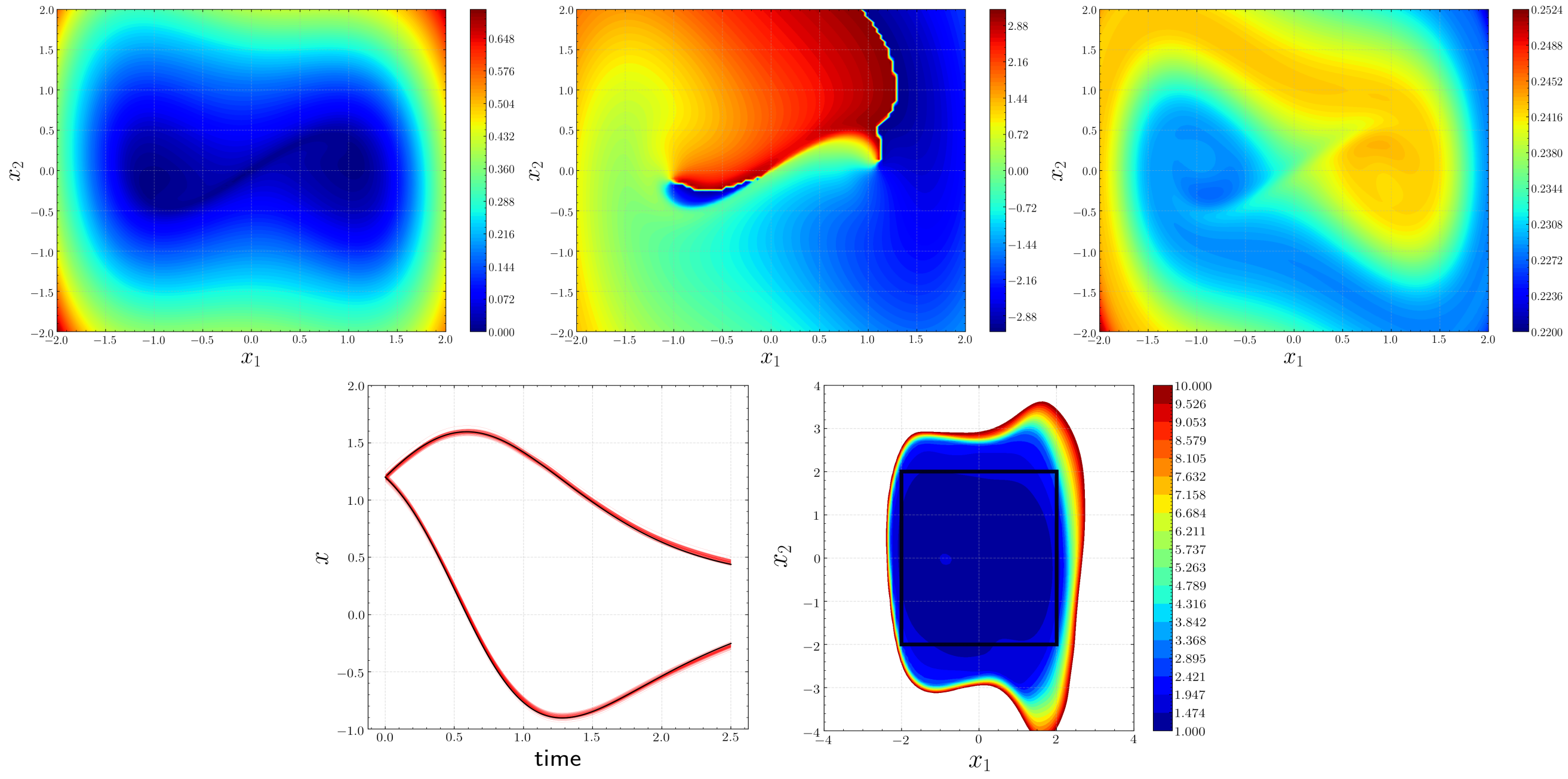}
    \caption{2D unforced Duffing oscillator system case: mean of Koopman observables and predicted trajectory from the differential form model. Top left: mean of the amplitude of Koopman eigenfunction associated with $\lambda_{1,2}$. Top middle: mean of the phase angle of Koopman eigenfunction associated with $\lambda_{1,2}$. Top right: mean of the amplitude of Koopman eigenfunction associated with $\lambda_3$. Bottom left: Monte Carlo sampling of the predicted trajectory with initial condition $\mathbf{x}_0 = [1.2,1.2]$. Bottom right: contour of the normalized standard deviation over $x_1, x_2 \in [-4,4]$ where the black square represents the boundary of training data.}
    \label{fig:2d_duffing_one_figure}
\end{figure}

To assess the uncertainty in the Koopman eigenfunctions, we draw 100 samples from the variational posterior of $\mathbf{W}_{\mathbf{\Phi}}$ and $\mathbf{K}$. 
First, the mean of the non-trivial Koopman eigenvalues are $\lambda_{1,2} = -0.535 \pm 0.750i$ and the mean of the third eigenvalue is $\lambda_{3}=-2 \times 10^{-5}$.
The mean of the module and phase angle of the Koopman eigenfunctions is shown in \cref{fig:2d_duffing_one_figure}. The results are similar to Ref.~\cite{otto2017linearly} in which a deterministic model is employed. The Koopman eigenfunction associated with $\lambda_3$ acts as an indicator of the basin of attraction.
Second, to better visualize the effect of uncertainty on unseen data, we normalize the standard deviation of the module of Koopman eigfunctions by the minimal standard deviation over $[-4,4]$.  \cref{fig:2d_duffing_one_figure} shows a uniformly distributed standard deviation within  $[-2,2]\times [-2,2]$ where sampling data is distributed, and a drastic increase outside that training region as expected. Note that the area where the normalized standard deviation is larger than ten is cropped. 

To obtain the posterior distribution of the evolution of the dynamics predicted by the Koopman operator, we arbitrarily choose an initial condition within the range of training data as $\mathbf{x}(0) = [1.2, 1.2]$. Note that we did not input the model with any trajectory data other than the state and the corresponding time derivative obtained from the governing equation. With Monte Carlo sampling, we obtain the distribution of the trajectory in \cref{fig:2d_duffing_one_figure}. Note that the uncertainty is quite small since it is within the training data region with enough data.

\subsection{Flow wake behind a cylinder}

We consider the velocity field of a two-dimensional laminar flow past a cylinder \textcolor{black}{in a transient regime} at a Reynolds number $Re_D =  U_{\infty} D / \nu=100$ \textcolor{black}{, which is above 47, the critical Reynolds number associated with Hopf bifurcation.}   $U_{\infty}$ is the freestream velocity, $D$ is the cylinder diameter, and $\nu$ is the kinematic viscosity. \textcolor{black}{The transient regime is rather difficult due to the high nonlinearity of post-Hopf-bifurcation dynamics between unstable equilibrium and the limit cycle~\cite{chen2012variants}.} Data is generated by solving the two dimensional incompressible Navier--Stokes equations using OpenFOAM~\cite{jasak2007openfoam}. \textcolor{black}{Grid convergence was verified using a sequence of successively refined meshes}. 

\textcolor{black}{The initial condition is a uniform flow with the freestream velocity superimposed with standard Gaussian random noise, and with pressure initialized with Gaussian random noise. Although this initial condition is a rough approximation to the development of true instabilities  from equilibrium~\cite{chen2012variants}, the flow is observed to rapidly converge to a quasi-steady solution after a few steps, and then starts to oscillate and form a long separation bubble with two counter-rotating vortices.} The first 50 POD snapshots of velocity are considered with the kinetic energy captured upto 99\%. We sample 1245 snapshots of data on the trajectory starting from the unstable equilibrium point to the vortex shedding attractor with $\Delta t = 0.1 t_{ref} = 0.1 D/U_{\infty}$ where the characteristic advection time scale $t_{ref}=D/U_{\infty}= 2\textrm{ sec}$.  The first 600 snapshots are  considered as training data, i.e., $0 \le t \le 60 t_{ref}$. 

To further analyze the robustness of the model to noisy data,  Gaussian white noise is added~\footnote{Note that  noise was not added to the original flow field since taking the dominant POD modes on the flow field would contribute to de-noising.} to the temporal data of POD coefficients by considering a fixed signal-to-noise ratio as 5\%, 10\%, 20\%, 30\%, for each component at each time instance. The model performance is evaluated by predicting the entire trajectory with the (noisy) initial condition given, including the remaining 645 snapshots. 
We consider the recurrent form model described in \cref{sec:recurrent} together with SVD-DMD described in \cref{sec:svd_dmd} with the corresponding hyperparameters given in 
\cref{tab:cylinder}. Note that we consider 20 intrinsic modes, i.e., at most 10 distinct frequencies can be captured, which are empirically chosen. Also, the finite-horizon window length corresponding to $T=100$ is \textcolor{black}{$10t_{ref}$}, which is much less than the time required for the system to transit from unstable equilibrium to the attractor.   
\begin{table}[htbp]
{   
    \footnotesize
    \caption{Hyperparameters of recurrent form model for flow past cylinder at $Re_{D}=100$.}
    \label{tab:cylinder}
    \begin{center}
    \begin{tabular}{|c|c|c|c|c|c|} \hline
        \bf layer structure & \bf optimizer  & \bf learning rate & \bf total epoch &  \bf batch size & \bf $T$ \\  \hline
          \textcolor{black}{50-100-50-20-20-20-50-100-50} &  Adam & 1e-3 & 20000 & 64 & 100 \\ \hline
    \end{tabular}
    \end{center}
}
\end{table}

\begin{figure}[htbp]
\centering
\includegraphics[width=0.7\linewidth]{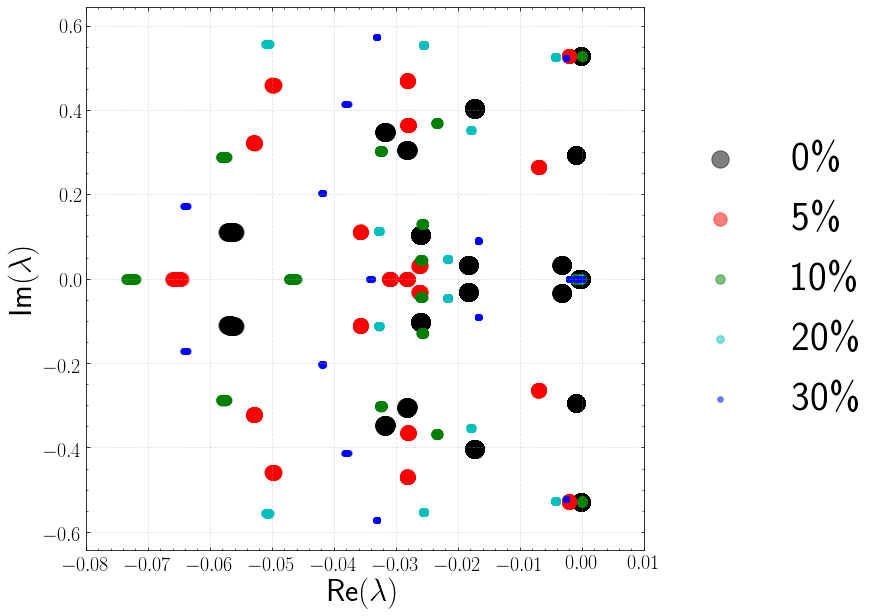}
\caption{Distribution of continuous-time Koopman eigenvalues for fixed signal-to-noise ratio from 5\%, 10\%, 20\%, 30\%.}
\label{fig:cyd_koopman}
\end{figure}

The continuous-time Koopman eigenvalue distribution of 20 modes for the training data with five different signal-to-noise ratios is shown in \cref{fig:cyd_koopman}. First, all Koopman eigenvalues are stable according to the stable parameterization of $\mathbf{K}$ in \cref{eq:stable_structure}. Second, when noise is added, the eigenvalues are seen to deviate except the one on the imaginary axis with $\lambda = \pm 0.528j$, which corresponds to the dominant vortex shedding frequency on the limiting cycle with $St = \lambda D/(2\pi U_{\infty}) = 0.168$. \textcolor{black}{This exercise verifies the robustness of the present approach to  Gaussian noise.}

\begin{figure}[htbp]
\centering
\includegraphics[width=\linewidth]{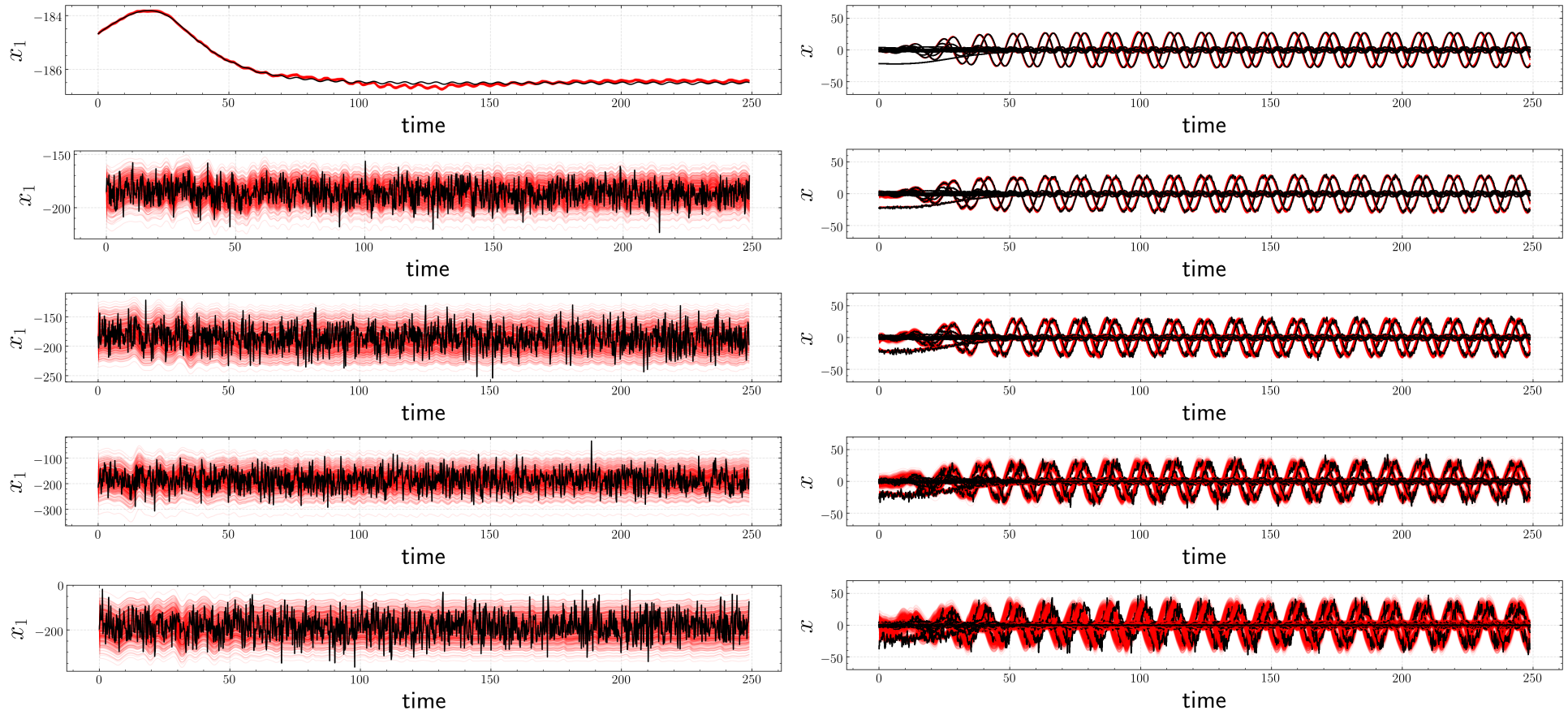}
\caption{Comparison between Monte Carlo sampled distribution of predicted trajectory (Red) and the noisy training data (Black) for signal-to-noise ratio from top to bottom as 0\% 5\%, 10\%, 20\%, 30\%. Left: first POD coefficient. Right: rest 49 POD coefficients.}
\label{fig:pred_pod_traj}
\end{figure}

\begin{figure}[htbp]
\centering
\includegraphics[width=0.8\linewidth]{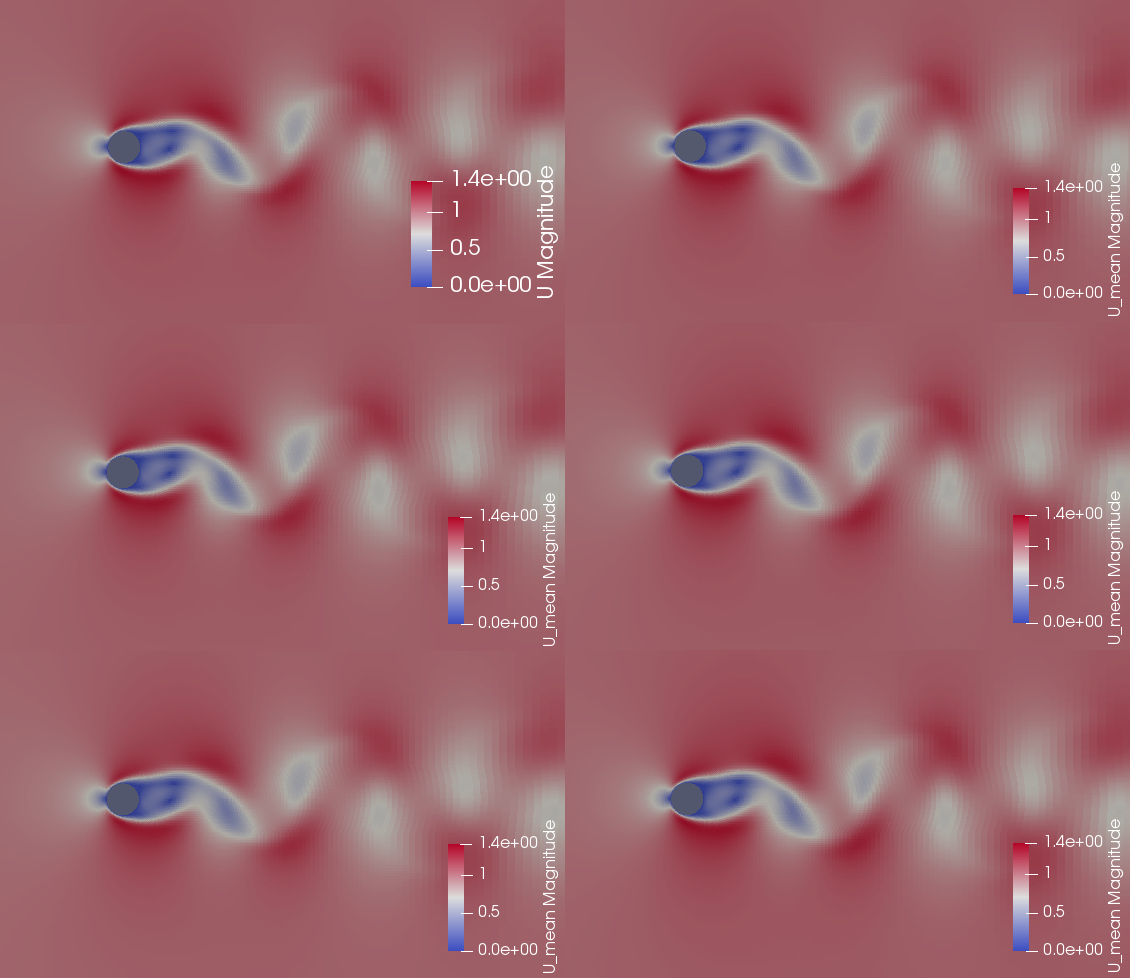}
\caption{Contour of velocity magnitude at $t=100t_{ref}$. Top left: flowfield projected from the true dominant 50 POD coefficients. Top right: flowfield projected from the 50 POD coefficients of the mean of posterior distribution with clean data. Middle left and right, bottom left and right correspond to noisy training data with signal-to-noise ratio 5\%, 10\%, 20\%, 30\% respectively.}
\label{fig:pred_mean_traj}
\end{figure}

\begin{figure}[htbp]
\centering
\includegraphics[width=0.8\linewidth]{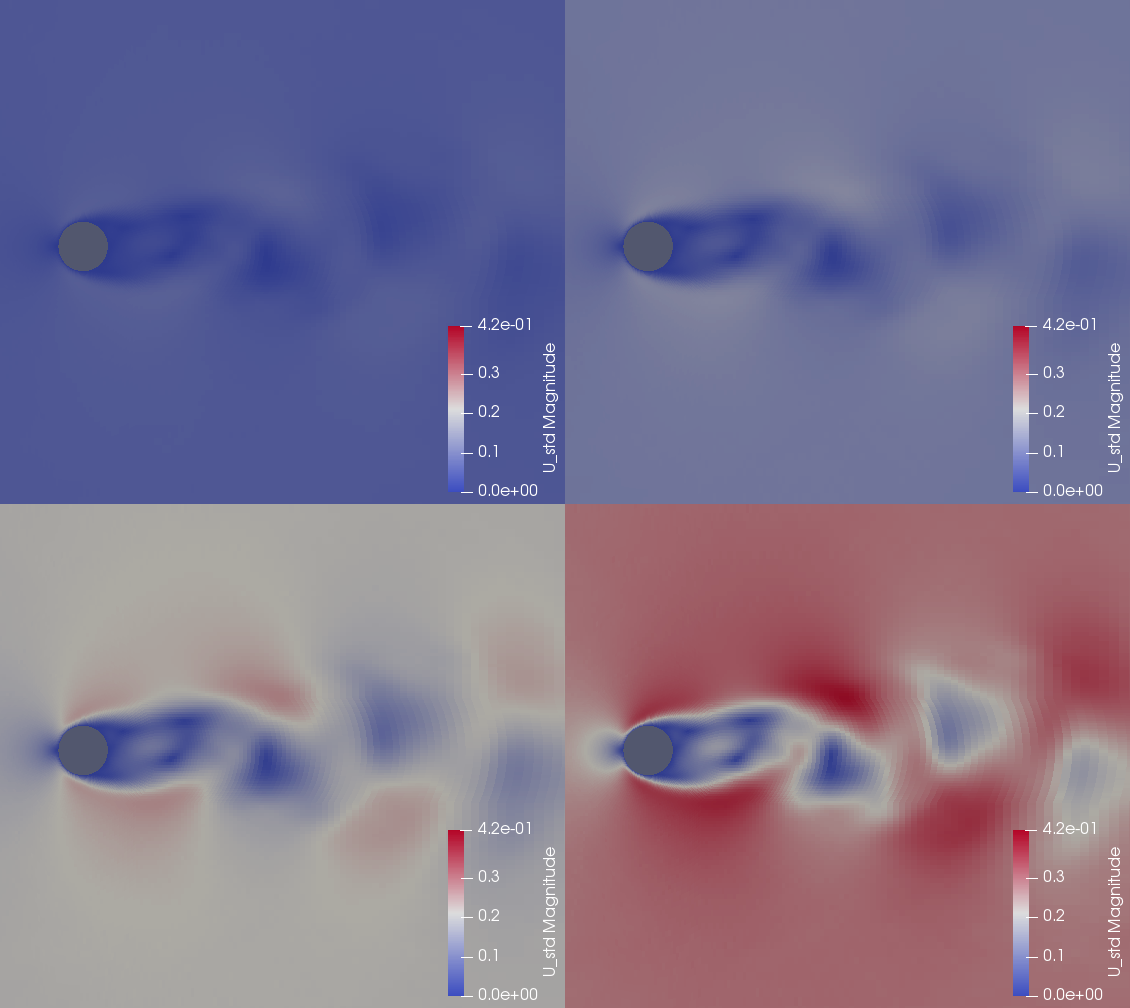}
\caption{Contour of standard deviation of velocity magnitude at $t=100t_{ref}$. Top left and right, bottom left and right correspond to noisy training data with signal-to-noise ratio 5\%, 10\%, 20\%, 30\% respectively.}
\label{fig:pred_std_traj}
\end{figure}

Recall that we have a posterior distribution of the predicted trajectory of POD components for  with five different noise ratios. Monte Carlo sampling of the posterior distribution is shown for the noisy training data in \cref{fig:pred_pod_traj}.
Clearly, uncertainty from the data due to the Gaussian white noise is well characterized by the ensemble of Monte Carlo sampling on the distribution. To further analyze the effect of the difference between the mean of the posterior distribution of the prediction and the ground truth clean trajectory on the flow field, we show the (projected) mean and standard deviation of the predicted POD coefficients at $t=100 t_{ref}$, in \cref{fig:pred_mean_traj} and \cref{fig:pred_std_traj}. As seen in \cref{fig:pred_mean_traj}, there is hardly any difference between the mean of the posterior distribution and the ground truth. The contour of standard deviation projected onto the flow field shows a similar pattern to the vortex shedding, and the standard deviation near the wake region is relatively small compared to other domains in the flow field. 

\section{Conclusions}
\label{sec:conclusion}

A probabilistic, stabilized deep learning framework was presented towards the end of extracting the Koopman decomposition for continuous nonlinear dynamical systems. We formulated the deep learning problem from a measure-theoretic perspective with a clear layout of all the assumptions. Two different forms: differential and recurrent, suitable for different situations were proposed and discussed. A  parameterization of Koopman operator was proposed, with guaranteed stability.  Further, a novel deep learning architecture was devised such that the SVD-DMD~\cite{schmid2010dynamic} is naturally embedded. Finally,  mean-field variational inference was used to quantify uncertainty in the modeling. To evaluate the posterior distribution, Monte Carlo sampling procedures corresponding to different forms were derived. Finally, the model was evaluated on three continuous nonlinear dynamical systems ranging from toy polynomial problems to an unstable wake flow behind a cylinder, from three different aspects: the uncertainty with respect to the density of data in the domain, unseen data, and noise in the data. The results show that the proposed model is able to capture the uncertainty in all the above cases and is robust to noise.  

\appendix
\section{Real part of eigenvalues of a negative semi-definite real square matrix}
\label{apdx}

\begin{definition}
\label{def:negative_semi}
An $n \times n$ real matrix (and possibly non-symmetric) $\mathbf{A}$ is called negative semi-definite, if $\mathbf{x}^\top \mathbf{A} \mathbf{x} \le 0$, for all non-zero vectors $\mathbf{x} \in \mathbb{R}^{n\times 1}$.  
\end{definition}

\begin{lemma}
\label{lem:lemma_negative_semi}
For a negative semi-definite real square matrix, the real part of all of its eigenvalues is non-positive.  
\end{lemma}

\begin{proof}
Consider the general negative semi-definite matrix defined in \cref{def:negative_semi} as $\mathbf{A} \in \mathbb{R}^{n \times n}$, such that $\mathbf{x}^\top \mathbf{A} \mathbf{x} \le 0$ for any non-zero vector $\mathbf{x} \in \mathbb{R}^{n\times 1}$. Without loss of generality, denote its eigenvalues as $\lambda = \alpha + j\beta$, where $\alpha, \beta \in \mathbb{R}$ and corresponding eigenvectors as $\mathbf{v} = \mathbf{v}_r + j\mathbf{v}_i$ where $\mathbf{v}_r, \mathbf{v}_i \in \mathbb{R}^n$. 

Then we have $0 = (\mathbf{A}- \lambda) \mathbf{v} = (\mathbf{A} - \alpha - j\beta)(\mathbf{v}_r + j\mathbf{v}_i)$, which further leads to $(\mathbf{A} - \alpha)\mathbf{v}_r = -\beta \mathbf{v}_i$ and $(\mathbf{A} - \alpha)\mathbf{v}_i = \beta\mathbf{v}_r$. Then we have $\mathbf{v}_r^\top (\mathbf{A} - \alpha)\mathbf{v}_r = -\beta \mathbf{v}_i^\top  \mathbf{v}_r$, and $\mathbf{v}_i^\top (\mathbf{A} - \alpha)\mathbf{v}_i =\mathbf{v}_i^\top \beta\mathbf{v}_r$. Thus adding the two equations, we have $\mathbf{v}_r^\top (\mathbf{A} - \alpha)\mathbf{v}_r + \mathbf{v}_i^\top (\mathbf{A} - \alpha)\mathbf{v}_i=0$, which implies $\alpha = (\mathbf{v}_r^\top \mathbf{A} \mathbf{v}_r + \mathbf{v}^\top_i \mathbf{A} \mathbf{v}_i )/(\mathbf{v}_r^\top \mathbf{v}_r + \mathbf{v}_i^\top \mathbf{v}_i)$. 

Using the definition of negative semi-definite matrices, we have $\mathbf{v}_r^\top \mathbf{A} \mathbf{v}_r \le 0 $ and $\mathbf{v}_i^\top \mathbf{A} \mathbf{v}_i \le 0$, even if $\mathbf{v}_r$ or $\mathbf{v}_i = \mathbf{0}$. Since there is at least one non-zero vector between $\mathbf{v}_r$ and $\mathbf{v}_i$, one can safely arrive at $\alpha \le 0$ that the real part of any eigenvalue of $\mathbf{A}$ is non-positive.
\end{proof}

\section*{Acknowledgements}
This work was supported by DARPA under the grant titled {\em Physics Inspired Learning and Learning the Order and Structure Of Physics,}.
Computing resources were provided by the NSF via grant 1531752 MRI: Acquisition of Conflux, A Novel Platform for Data-Driven Computational Physics.

\bibliographystyle{siamplain}
\bibliography{references}

\end{document}